\documentclass[12pt,letterpaper]{article}
\usepackage{setspace}
\usepackage[utf8]{inputenc}
\usepackage{graphicx}
\usepackage{amsmath}
\usepackage{amssymb,amsthm}
\usepackage{geometry}
\usepackage{cases}
\usepackage{titlesec}
\graphicspath{ {images/} }
\newtheorem{theorem}{Theorem}
\newtheorem{lemma}{Lemma}
\newtheorem{cor}{Corollary}
\newtheorem{prop}{Proposition}
\newtheorem{de}{Definition}
\newtheorem{rk}{Remark}

\DeclareMathOperator*{\supess}{ess\,sup}
\DeclareMathOperator*{\supp}{supp}
\titleformat{\section}[hang]
{\Large\bfseries}
{\thesection.}{0.5em}{}
\begin{document}
\def\R{\mathbb{R}}
\def\H{\mathbb{H}}
\def\C{\mathbb{C}}
\def\a{\alpha}
\def\e{\epsilon}
\def\ds{\displaystyle}
\def\ni{\noindent}
\def\avint{\mathop{\,\rlap{\bf--}\!\!\int}\nolimits}
\def\ba{\begin{aligned}}
\def\ea{\end{aligned}}
\def\ol{\overline}
\def\wt{\widetilde}
\def\b{\beta}\def\ab{{\overline\a}}
\def\avgdj{\avint_{D_{j+1}\setminus D_{j}}}
 \def\avgdJ{\avint_{D_{J+1}\setminus D_{J}}}
 \def\avgdJo{\avint_{D_{J_1+1}\setminus D_{J_1}}}
 \def\dJ{\int_{D_{J+1}\setminus D_{J}}}
 \def\dj{\int_{D_{j+1}\setminus D_{j}}}
 \def\wphi{\wt{\phi}_{\e,r}}
 \def\mc{\mathcal}
 \def\wh{\widehat}
\def\be{\begin{equation}}\def\ee{\end{equation}}\def\nt{\notag}\def\bc{\begin{cases}}\def\ec{\end{cases}}\def\ba{\begin{aligned}}
\def\ea{\end{aligned}}
\def\QEDopen{{\setlength{\fboxsep}{0pt}\setlength{\fboxrule}{0.2pt}\fbox{\rule[0pt]{0pt}{1.3ex}\rule[0pt]{1.3ex}{0pt}}}}
\centerline{\bf \large{Adams inequalities with exact growth condition }}
\centerline{\bf\large{for Riesz-like potentials  on $\R^n$}}
\vskip 0.1in
\centerline{\bf{Liuyu Qin}}

%
\begin{abstract}
We derive sharp Adams inequalities with exact growth condition for the Riesz potential and for more general Riesz-like potentials on $\R^n$. We also obtain Moser-Trudinger inequalities with exact growth condition for the fractional Laplacian, and for general homogeneous elliptic differential operators with constant coefficients.
\end{abstract}
\numberwithin{equation}{section}



\def\R{\mathbb{R}}
\def\H{\mathbb{H}}
\def\C{\mathbb{C}}
\def\a{\alpha}
\def\e{\epsilon}
\def\ds{\displaystyle}
\def\ni{\noindent}
\def\avint{\mathop{\,\rlap{\bf--}\!\!\int}\nolimits}
\def\ba{\begin{aligned}}
\def\ea{\end{aligned}}
\def\ol{\overline}
\def\wt{\widetilde}
\def\b{\beta}\def\ab{{\overline\a}}
\def\avgdj{\avint_{D_{j+1}\setminus D_{j}}}
 \def\avgdJ{\avint_{D_{J+1}\setminus D_{J}}}
 \def\avgdJo{\avint_{D_{J_1+1}\setminus D_{J_1}}}
 \def\dJ{\int_{D_{J+1}\setminus D_{J}}}
 \def\dj{\int_{D_{j+1}\setminus D_{j}}}
 \def\wphi{\wt{\phi}_{\e,r}}
 \def\mc{\mathcal}
 \def\wh{\widehat}
\def\be{\begin{equation}}\def\ee{\end{equation}}\def\nt{\notag}\def\bc{\begin{cases}}\def\ec{\end{cases}}\def\ba{\begin{aligned}}
\def\ea{\end{aligned}}
\def\QEDopen{{\setlength{\fboxsep}{0pt}\setlength{\fboxrule}{0.2pt}\fbox{\rule[0pt]{0pt}{1.3ex}\rule[0pt]{1.3ex}{0pt}}}}
\section{Introduction and main results}
\par
\ \ \ 

The  Moser-Trudinger inequality with exact growth condition on $\R^n$ takes the form

\be \int_{\R^n} \frac{\exp_{\lceil\frac{n}{\a}-2\rceil}\left[\b(\a,n)|u|^{\frac{n}{n-\a}}\right]}{1+|u|^{\frac{n}{n-\a}}}dx\le C||u||^{n/\a}_{n/\a}\ \ \quad  \text{for all }\ u\in W^{\a,\frac{n}{\a}}(\R^n),\ \ ||\nabla^\a u||_{n/\a}\le 1\label{in17}\ee
where $\lceil x\rceil$ denotes the ceiling of $x$, i.e. the smallest integer greater than or equal $x$, and where $\exp_N$ is the regularized exponential, that is
\be \exp_N(t)=e^t-\sum_{k=0}^{N}\frac{t^k}{k!}\nt\ee
and where for $\alpha\in(0,n)$ an integer the higher order gradient $\nabla^\a$ is defined as 
\be \nabla^\a u=\bc (-\Delta)^{\frac \a 2}u\ \ \ &\text{if}\ \a\ \text{is even}\\ \nabla (-\Delta)^{\frac {\a-1} 2}u\ \ \ &\text{if}\  \a\ \text{is odd.}\ec\nt\ee
 
Such inequality was proved first  by Ibrahim, Masmoudi and Nakanishi  [IMN] for $n=2$ and $\a=1$, followed by Masmoudi and Sani  who dealt with the cases $n=4$, 
$\a=2$ in [MS1],    any $n$ and $\alpha=1$ in [MS2], and any $n$ and any integer  $\alpha$ in [MS3]. In [LT] Lu and Tang dealt with the case $\alpha=2$, for any $n$.
  In all these results the explicit sharp exponential constant $\beta(\a,n)$ (see \eqref{5c}) that was found  was the same as the sharp exponential constant in the classical Moser-Trudinger  inequality on bounded domains due to Adams [A1]:

 \be \int_{\Omega} \exp\left[\b(\a,n)|u|^{\frac{n}{n-\a}}\right]dx\le C|\Omega|\ \ \quad  \text{for all } u\in W_0^{\a,\frac{n}{\a}}(\Omega),\ \ \|\nabla^\a u\|_{n/\a}\le 1\label{Adams}\ee
where $|\Omega|$ denotes the Lebesgue measure of $\Omega$. Recall that the exponential constant is sharp in the sense  that it cannot be replaced by a larger constant.
\def\na{\frac n{\a}}

The main new result behind the proof of \eqref{in17} in [IMN], [MS1] ($\alpha=1$) is what the authors call ``optimal descending growth condition" (ODGC). In essence,  such result gives an optimally adjusted exponential growth of {\it radial} functions outside balls of radius $R$, given the $L^n$ norms of their gradients. In [MS2] and  [LT]   the same result is proven for radial functions under $L^{n/2}$ norm conditions on  their Laplacians, and in [MS3]  under Lorentz $L^{n/2,q}$ norm conditions on their Laplacians. The key initial step  that allowed the authors to only consider the radial Sobolev functions was the application of   well-known, powerful symmetrization inequalities, specifically the P\'olya-Szeg\"o and Talenti's inequalities.

There are a few other types of sharp Moser-Trudinger inequalities in the whole $\R^n$. The most common one states that for all $u\in W^{\a,\na}(\R^n)$ satisfying  the under the {\it Ruf condition}
\be \|u\|_{n/\a}^{n/\a}+\|\nabla^\a u\|_{n/\a}^{n/\a}\le 1\label{ruf}\ee
the following estimate holds
\be \int_{\R^n}\exp_{\lceil\frac{n}{\a}-2\rceil}\left[\b(\a,n)|u|^{\frac{n}{n-\a}}\right]dx\le C.\label{MT1}\ee

This result was first derived by Ruf, [R] for $\alpha=1$ in dimension $n=2$ and later extended to all dimensions by Li-Ruf [LR]. The general case was settled by Fontana-Morpurgo in [FM2], where the authors prove \eqref{MT1} under \eqref{ruf} for arbitrary $n$ and integer $\alpha$, but also for fractional powers of $\Delta$, and for homogeneous elliptic operators with constant coefficients. 

Under norm conditions weaker than \eqref{ruf} estimate \eqref{MT1} is in general false, but it becomes true if one lowers the exponential constant. For example under the condition 
\be \max\big\{\|u\|_{n/\a}, \,\|\nabla^\a u\|_{n/\a}\big\}\le 1\label{at}\ee
inequality \eqref{MT1} holds with exponential constant $\theta\beta(\a,n)$, for any $\theta\in (0,1)$. This result was originally derived for $\alpha=1$  by Cao [C] and Panda [P]  in dimension 2 and by Do \'O [D] in any dimension. Later Adachi-Tanaka [AT] re-proved the result and cast it in a dilation invariant form. In [FM2] the authors derived the general case as a corollary of \eqref{MT1} under \eqref{ruf},  and showed that under either \eqref{at} or under 

\be \|u\|_{n/\a}^{rn/\a}+\|\nabla^\a u\|_{n/\a}^{rn/\a}\le 1,\qquad r>1\label{at1}\ee
for any $\theta\in(0,1)$

\be \int_{\R^n}\exp_{\lceil\frac{n}{\a}-2\rceil}\left[\theta\b(\a,n)|u|^{\frac{n}{n-\a}}\right]dx\le C(1-\theta)^{-1+1/r},\label{MT2}\ee
where $r=\infty$ under \eqref{at}.

It is important to point out that the Masmoudi-Sani  result is the strongest one to date, in the sense that it directly implies \eqref{MT1} under the Ruf condition (see [MS1], [MS2], [MS3]).


Our initial goal was to derive  the sharp Adams inequality with exact growth condition for the Riesz potential
$$I_\a f(x)=\int_{\R^n} |x-y|^{\a-n} f(y)dy,$$
that is
\be\int_{\R^n}\frac{\exp_{\lceil\frac{n}{\a}-2\rceil}\left[\dfrac{1}{|B_1|}|I_\a f|^{\frac{n}{n-\alpha}}\right]}{1+|I_\a f|^{\frac{ n}{n-\alpha}}}dx \leq C||I_\a f||_{n/\alpha}^{n/\alpha},\qquad ||f||_{\frac{n}{\a}}\le 1, \label{in18}\ee
 where  $|B_1|$ is the volume of the unit ball of $\R^n$ and where the exponential constant is sharp. Note that the exponential constant $|B_1|^{-1}$ is the same constant as in the original inequality due to Adams [A1]:
 \be \int_{\Omega} \exp\left[\dfrac{1}{|B_1|}|I_\a f(x)|^{\frac{n}{n-\a}}\right]dx\le C|\Omega|\ \ \quad  \text{for all } f\in L^{\frac{n}{\a}}(\Omega),\ \ \|f\|_{n/\a}\le 1\label{Adams1}.\ee
Clearly \eqref{in18} implies \eqref{in17}, in the same way that \eqref{Adams1} implies \eqref{Adams} due to the fact that $I_\a$ is the inverse of $(-\Delta)^{\a/2}$ on smooth, compactly supported functions.

In this paper  we prove that \eqref{in18} is true, and not only for the Riesz kernel but for a subclass of the Riesz-like  kernels introduced by Fontana-Morpurgo, which we call \emph{strictly Riesz-like kernels.}\par
\bigskip
To describe our result let us recall the definition given in [FM2]:
\begin{de}\label{rzlike} A measurable function $K\ :\ \R^n\setminus \{0\}\to \R$ is a Riesz-like kernel of order $\a\in(0,n)$ if it satisfies the following properties:
 \be K(x)=g(x^*)|x|^{\alpha-n}+O(|x|^{\alpha-n+\delta_1}) \ ,\qquad  x^*=\frac{x}{|x|},\qquad 0<|x|\le B\label{A1}\tag{A1}\ee 
 \be|K(x)|\leq H_1|x|^{\alpha-n}\label{A2}\tag{A2}\ee
 \be|K(z_1)-K(z_2)|\leq H_2|z_1-z_2|\max \{{|z_1|}^{\alpha-n-1},{|z_2|}^{\alpha-n-1}\},\ \ z_1, z_2\neq0 \label{A3}\tag{A3}\ee
 where $g\ :\ S^{n-1}\to\R$ is a measurable function and $\delta_1, H_1,H_2,B$ are positive constants.
\end{de}
If we add an additional condition \eqref{A4} as below, we have more restrictive control of the kernel $K$ when $|x|$ is large:
\begin{de}\label{riesz2} A measurable function $K\ :\ \R^n\setminus \{0\}\to \R$ is a strictly Riesz-like kernel of order $\a\in(0,n)$ if it is Riesz-like and satisfies the following property:
 \be |K(x)|\le |g(x^*)||x|^{\alpha-n}+O(|x|^{\alpha-n-\delta_2}) \ ,\ \ |x|>B \label{A4}\tag{A4}\ee
 where $g\ :\ S^{n-1}\to\R$ is a measurable function and $B, \delta_1,\delta_2, H_1,H_2$ are positive constants.
\end{de}
Here the ``big O'' notation in \eqref{A1} means that $|O(|x|^{\alpha-n+\delta_1})|\le C|x|^{\alpha-n+\delta_1}$ for all $x$ such that $0<|x|\le B$. And the same notation in \eqref{A4} means that $|O(|x|^{\alpha-n-\delta_2})|\le C|x|^{\alpha-n-\delta_2}$.  It is clear that \eqref{A3} implies that $g$ is Lipschitz. Also, \eqref{A1},\eqref{A3} and \eqref{A4} imply \eqref{A2}.

Clearly, any kernel of type $g(x^*)|x|^{\a-n}$ with $g$ Lipschitz on the sphere, provides an example of strictly Riesz-like kernel.
\par

For $m\in \mathbb{N}$, a kernel $K$ is called \emph{m-regular} if $K\in C^m(\R^n\setminus \{0\})$ and 
\be |D_x^h K(x)|\le C|x|^{\a-n-|h|},\ \ \ x\ne 0,\ |h|\le m\nt\ee
where $h=(h_1,...,h_n)$ is a multi-index with $|h|=h_1+...+h_n$ and where $D_x^h K$ denotes the $h$-th derivative of $K$ with respect to $x$. Clearly the Riesz kernel is $m$-regular for all $m$, and any $1$-regular $K$ satisfies condition \eqref{A3}. \par
  Let us denote $T$ the convolution operator with kernel $K$:
  \be Tf(x)=K\ast f(x)=\int_{\R^n} K(x-y)f(y)dy.\nt\ee
  For vector valued functions
$K=(K_1,...,K_m),\ f=(f_1,...,f_m)$
  we define $Tf$ in the same way with
  \be Kf=K_1f_1+...+K_mf_m,\ \ \ |f|=(f_1^2+...+f_m^2)^{1/2}.\nt\ee
  \par The results and proofs in this paper apply to both scalar and vector cases, so we will not distinguish between these two cases, except in the proof of sharpness (see Remark \ref{vectorsharp}).\par
 

\bigskip
The main result of this paper is the following:
 \begin{theorem}\label{T1}
 Let $0<\alpha <n$, and $K$ is strictly Riesz-like. There exists $C=C(n,\alpha,K)>0$ such that for all compactly supported $f$ with
 $$ ||f||_{\frac{n}{\alpha}}\leq 1$$ we have
 \be\int_{\R^n}\frac{\exp_{\lceil\frac{n}{\a}-2\rceil}\bigg[\dfrac{1}{A_g}|Tf|^{\frac{n}{n-\alpha}}\bigg]}{1+|Tf|^{\frac{ n}{n-\alpha}}}dx \leq C||Tf||_{n/\alpha}^{n/\alpha}, \label{1b}\ee
 where \be A_g=\dfrac 1 n \int_{S^{n-1}}|g(x^*)|^{\frac{n}{n-\alpha}}dx^*.\label{ag1}\ee 
If $K$ is $n$-regular, then the exponential constant $A_g^{-1}$ in  \eqref{1b} cannot be replaced by a larger number.
 Furthermore, if $K$ is $n$-regular, then \eqref{1b} cannot hold if the power $\frac{n}{n-\a}$ in the denominator is replaced by any smaller power.
\end{theorem}
Here $dx^*$ is the surface measure of the unit sphere $S^{n-1}$, induced by the Lebesgue measure.

\smallskip \newpage
As pointed out in [FM2] Adams type  estimates involving an integral of the regularized exponential  over the whole space, have  equivalent ``local" formulations in terms of the standard exponential. Via the exponential regularization lemma (Lemma A in section 3) estimate \eqref{1b} is equivalent to the following  local version
\be\int_{E}\frac{\exp\bigg[\dfrac{1}{A_g}|Tf|^{\frac{n}{n-\alpha}}\bigg]}{1+|Tf|^{\frac{ n}{n-\alpha}}}dx \leq C\big(|E|+||Tf||_{n/\alpha}^{n/\alpha}\big) \label{1c}\ee
valid for all measurable $E$ with finite measure, and under $ ||f||_{\frac{n}{\alpha}}\leq 1$.\par
We mention that inequality \eqref{1b} still holds if we have ``$\le$''  instead of ``$=$'' in condition \eqref{A1}, that is, if \eqref{A1} is replaced by
\be |K(x)|\le |g(x^*)||x|^{\alpha-n}+|O(|x|^{\alpha-n+\delta_1})| ,\qquad 0<|x|\le B.\label{A1'}\tag{A1'}\ee 
But in order to have sharpness in the exponential constant $A_g^{-1}$ we need to assume condition \eqref{A1}.\par

We point out that for Theorem \ref{T1} to hold it is not enough to assume that $K$ be only Riesz-like. It is relatively easy to find an example of a Riesz-like kernel such that the inequality in Theorem \ref{T1} cannot hold, but the one in [FM2, Theorem 5] holds. In section \ref{sharpness} remark \ref{example}, we will address this example, which indicates that it is necessary for us to strengthen our assumption for large $|x|$, i.e. \eqref{A4}, so that $K$ has same behavior near  the origin  and at infinity.  \par

One of the main difficulties we had to overcome toward a proof of Theorem 1, even for the Riesz potential as in \eqref{in18}, was to find  a suitable replacement of the optimal growth condition result  for the potential $Tf$, under norm conditions on $f$. Clearly, in this context no tools such as the P\'olya-Szeg\"o or Talenti's
inequalities are available, which makes an initial reduction to radial functions impossible, even in the case of the Riesz potential.  The way we bypass this problem is by carefully splitting the function $f$ and by making use of an improved O'Neil inequality. Very loosely speaking, we will consider a suitable 1-parameter family of sets $F_\tau$ depending on $f$, and  with measure $\tau$, and we will split $f$ as $f=f_\tau+f_\tau'$, with $f_\tau=f\chi_{F_\tau}^{}.$ By use of an improved O'Neil inequality we will prove an estimate of type 
\be (Tf)^*(t)\le U f_\tau(t)+ U' f_\tau' (\tau), \;\;0<t\le \tau \label{estimate}\ee
 where  $(Tf)^*$ is the symmetric decreasing rearrangement of $Tf$,  and where  $U, U'$ are two suitable, real-valued (nonlinear) functionals stemming from the O'Neil inequality (see estimate \eqref{e1}).  The first term in \eqref{estimate} is handled by an Adams inequality for sets of finite measure due to Fontana-Morpurgo (see Theorem A). The challenging part  is the proof of an optimal descending growth condition for the function $U'f_\tau'(\tau)$ (see Proposition 1).  In [IMN], [MS1], [MS2], [LT] a version  ODGC was first proved for sequences, followed by a suitable discretization of radial Sobolev functions. We will also make use of the discrete ODGC for sequences (See Lemma \ref{lo}, Section 4), however the discretization of $U'f'_\tau(\tau)$ turns out to be rather involved (see Proposition~1 and its proof, given in Section~5).


\bigskip

\par As a consequence  of Theorem \ref{T1} we derive the following general  Adachi-Tanaka type inequality: 
\begin{cor}\label{T3}
If  $K$ is a strictly Riesz-like kernel, then for any $\theta\in(0,1)$ there exists $C$ independent of $\theta$ such that for all compactly supported $f$ with  
\be ||f||_{n/\a}\le 1,
\label{3a}
\ee
we have 
\be \int_{\R^n}\exp_{\lceil\frac{n}{\alpha}-2\rceil}\bigg[\frac{\theta}{A_g}|Tf(x)|^{\frac{n}{n-\alpha}}\bigg]dx\leq \frac{C}{1-\theta}||Tf||_{n/\a}^{n/\a} \label{3c}\ee
where $A_g$ is the same as in \eqref{ag1}. If $K$ is $n$-regular and $K\notin L^{\frac{n}{n-\alpha}}(|x|\geq 1)$ then inequality \eqref{3c} is sharp, in the sense that the exponential integrals cannot be uniformly bounded if $\theta=1$.
\end{cor}

Estimate \eqref{3c} improves the one obtained in [FM2, Theorem 6], which does not have $\|Tf\|_{n/\a}$ on the right hand side, and which has $(1-\theta)^{-1}$ only in the case $K$ homogeneous. 
\def\na{\sum_{k=0}^{2N-1}}
\def\ia{\sum_{k=2N}^{\infty}}
\def\aa{\sum_{k=0}^{\infty}}

At the level of Moser-Trudinger inequalities, Theorem 1 implies almost immediately the Masmoudi-Sani result \eqref{in17}, for integer powers $\alpha$. 
Similarly, as a  consequence of Theorem 1, we will obtain  a Moser-Trudinger inequality with exact growth condition for the fractional Laplacian $(-\Delta)^{\a/2}$, for any $\a\in (0,n)$ , and also for general homogeneous elliptic operators.\par
To describe such result, recall that the Sobolev space $ W^{\a,p}({\R^n})$ consists of all locally summable functions $u\ :\ {\R^n}\to \R$ such that for each multiindex $h$ with $|h|\le \a$, the $h$-th weak partial derivative of $u$ exists and belongs to $L^p({\R^n}).$ For non integer $\a$, the space $W^{\a,\frac{n}{\a}}(\R^n)$ will denote  the Bessel potential space \be W^{\a,p}(\R^n)=\{u\in \mc{S}'\ :\ (I-\Delta)^{\a/2}u\in L^p(\R^n)\}=\{G_\a\ast f,\ f\in L^p(\R^n)\},\label{r2}\ee
where $G_\a$ is the kernel of the Bessel potential $(I-\Delta)^{-\a/2}$ and its Fourier transform is $(1+4\pi^2|\xi|^2)^{-\a/2}$.

We also recall that a homogeneous elliptic differential operator of even order $\a<n$ with real constant coefficients has the form
\be Pu=\sum_{|k|=\a}a_kD^ku\label{5e}\ee
for $u\in C_c^\infty(\R^n)$, with symbol
\be p_\a(\xi)=P(2\pi i \xi)=(2\pi)^\a(-1)^{\a/2}\sum_{|k|=\a}a_k\xi^k,\ \ |p_\a(\xi)|\geq c_0|\xi|^\a,\ \ \xi\in\R^n\nt\ee
for some $c_0>0$. The fundamental solution of $P$ is given by a convolution operator with kernel $g_P:$
\be g_P(x)=\int_{\R^n} \frac{e^{-2\pi i x\cdot \xi}}{p_\a(\xi)}d\xi \label{5f}\ee
in the sense of distributions. Since $p_\a$ is homogeneous of order $\a$, the kernel $g_P$ is also homogeneous with order $\a-n$. 
\begin{theorem}\label{T5}
For $0<\alpha<n, $ let $P$ be either $(-\Delta)^{\frac{\alpha}{2}}$, $\nabla(-\Delta)^{\frac{\alpha-1}{2}}$ for $\alpha$ odd, or a homogeneous elliptic operator of even order $\a<n$ with constant coefficients. Then there exists $C=C(\a,n,P)>0$ such that for every $u\in W^{\a,\frac{n}{\a}}(\R^n)$ with \be||Pu||_{n/\a}^{n/\a}\leq 1\label{5a}\ee
we have
\be \int_{\R^n}\frac{\exp_{\lceil\frac{n}{\alpha}-2\rceil}\big[\gamma(P)|u(x)|^{\frac{n}{n-\alpha}}\big]}{1+|u(x)|^{\frac{n}{n-\a}}}dx\leq C||u||_{n/\a}^{n/\a} \label{5b}\ee
where 
\be \gamma(P)=
\begin{cases}
\dfrac{c_\a^{-\frac{n}{n-\a}}}{|B_1|},\ \ \ &\text{if} \ P=(-\Delta)^{\frac{\alpha}{2}}\\
\dfrac{((n-\a-1)c_{\a+1})^{-\frac{n}{n-\a}}}{|B_1|},\ \ &\text{if} \ P=\nabla(-\Delta)^{\frac{\alpha-1}{2}}\ and \ \a\ odd,
\end{cases}
\label{5c}
\ee 
 with \be c_\a=\frac{\Gamma(\frac{n-\a}{2})}{2^\a \pi^{n/2}\Gamma(\frac{\a}{2})}\label{calpha}\ee
and where \be \gamma(P)=\frac{n}{\int_{S^{n-1}}|g^{}_P(x^*)|^{\frac{n}{n-\alpha}}dx^*}\nt\ee
if $P$ elliptic and $\a$ even. The exponential constant $\gamma(P)$ is sharp. Moreover, the above inequality \eqref{5b} cannot hold if the power $\frac{n}{n-\a}$ in the denominator is replaced by any smaller power.
\end{theorem}
 As an immediate consequences of Theorem \ref{T5}, we have the following Corollary:
 \begin{cor}\label{corT5}
 Let \  $\Omega$ be a bounded and open set in $\R^n$, $0<\a<n$ an integer. There exists $C>0$ such that for all $u\in W_0^{\a,\frac{n}{\a}}(\Omega)$ with $||\nabla^\a u||_{\frac{n}{\a}}\le 1$ we have
 \be  \int_{\Omega}\frac{\exp_{\lceil\frac{n}{\alpha}-2\rceil}\big[\gamma(P)|u(x)|^{\frac{n}{n-\alpha}}\big]}{1+|u(x)|^{\frac{n}{n-\a}}}dx\leq C||u||_{n/\a}^{n/\a}. \label{cT5a}\ee
 The exponential constant $\gamma(P)$ is sharp. Furthermore, the above inequality \eqref{cT5a} cannot hold if the power $\frac{n}{n-\a}$ in the denominator is replaced by any smaller power.
 \end{cor} 
 Although the proof of \eqref{cT5a} uses Adams inequality on $\Omega$, it is still not an easy direct consequence from Adams [A1]. We also mention that the above inequality \eqref{cT5a} is different from the inequalities in [A1] because of the norm of $u$ in the RHS.  For example, take $\a=1,\ n=2$, then by Corollary \ref{corT5} we have 
 \be  \int_{\Omega}\frac{e^{4\pi u^2}-1}{1+u^2}dx\leq C||u||_2^2 \label{cT5b}\ee
 and [A1] gave 
 \be  \int_{\Omega}{e^{4\pi u^2}}dx\leq C|\Omega| \label{cT5c}\ee
 So for fixed $\Omega$, we can see that, if $||u||_2^2$ becomes very small, then so is the LHS of \eqref{cT5b}, but this point may not be reflected by the second inequality \eqref{cT5c}.\par
 
\def\L{{\cal L}}
As we mentioned earlier, Riesz-like kernels were introduced in   [FM2], where  the authors proved, among other things,  that if  $K$ is a Riesz-like kernel, then under the Ruf condition  \be ||f||^{n/\alpha}_{n/\alpha}+||Tf||^{n/\alpha}_{n/\alpha}\leq 1\label{4a}\ee  the following Adams  inequality holds:
\be \int_{\R^n}\exp_{\lceil{\frac{n}{\alpha}-2\rceil}}\bigg[\frac{1}{A_g}|Tf(x)|^{\frac{n}{n-\alpha}}\bigg]dx\leq C \label{4c}\ee where $A_g$ is as in \eqref{ag1}.
 and where the exponential constant  $A_g^{-1}$ is  sharp  if the kernel is $n-$regular.
In section \ref{rufinq} we will prove that such result is implied by Theorem~\ref{T1} if $K$ is strictly Riesz-like.

\section{Improved O’Neil Lemma, O’Neil functional and Adams inequality}\label{ON}
\ \ \ Suppose that $(M,\mu)$ and $(N,\nu)$ are $\sigma$-finite measure spaces. Given a measurable function $f\ :\ M\rightarrow \ [-\infty,\infty]$ its distribution function will be denoted by
\be m_f(s)=\mu(\{x\in M:\ |f(x)|>s\}),\ \ \ s\geq 0.\nt\ee  
Assume that the distribution function of $f$ is finite for $s> 0$.

The decreasing rearrangement of $f$ will be denoted by
\be f^*(t)=\inf\{s\geq0:\ m_f(s)\leq t\},\ \ \ t>0\nt\ee
and we define \be f^{**}(t)=\frac{1}{t}\int_0^tf^*(u)du,\ \ \ t>0\nt\ee
which is sometimes called the \emph{maximal function} of $f^*$.
\par
Given a $\nu\times\mu$-measurable function $k: N\times M\rightarrow [-\infty,\infty]$, assume that the level sets of $k(x,\cdot)$ and $k(\cdot,y)$ have finite measure for all $x\in N$ and all $y\in M$. Let 
$k_1^*(x,t)$ and $k_2^*(y,t)$ be the decreasing rearrangement of $k(x,y)$ with respect to the variable $y$ (resp. $x$) for fixed $x$ (resp. $y$), and define
\be
\ba
&k_1^*(t)=\supess_{x\in N}k_1^*(x,t)\cr
&k_2^*(t)=\supess_{y\in M}k_2^*(y,t).
\ea
\nt\ee

Lastly, let $T$ be an integral operator defined as
\be Tf(x)=\int_M k(x,y)f(y) d\mu(y). \ee
One of the main tools used in the proof is the following slightly more general version of O'Neil lemma. 
\begin{lemma}\label{l}\textup{\textbf{(Improved O'Neil lemma)}}\\
Let  $k: N\times M\rightarrow [-\infty,\infty]$ be measurable, and 
\be k_1^*(t)\leq Dt^{-\frac{1}{\beta}},\ \ \ \ k_2^*(t)\leq Bt^{-\frac{1}{\beta}},\ \ \ t>0\label{l1}\ee
with $\beta>1$. Let $f:N\times M\to\R$ be a measurable function on $N\times M$. For each $x\in N$, let $f_x: M\to\R$ be defined as $f_x(y)=f(x,y)$ on $M$. Suppose there is a measurable function $\ol{f}:M\rightarrow [0, \infty]$ , $\ol{f}\in L^1(M)$ such that for $\nu-$a.e. $x\in N$
\be |f_x(y)|\leq \ol{f}(y), \ \ \mu - a.e.\  y\in M.\label{l2}\ee
Let 
\be T'f(x)=Tf_x(x)=\int_M k(x,y)f_x(y) d\mu(y), \label{l3}\ee
then $T'f(x)$ is well-defined and finite for  $\nu-$a.e. $x\in N$, and  there is a constant $C_0=C_0(D,B,\b)$ such that
\be (T'f)^{**}(t)\leq C_0\max\{\tau^{-\frac{1}{\beta}},t^{-\frac{1}{\beta}}\}\int_0^{\tau} \ol{f}^*(u)du +\supess_{x\in N} \int_\tau^\infty k_1^*(x,u)f_x^*(u)du,\ \ \ \forall t,\tau>0. \label{l4}\ee
\end{lemma}
Note that in this Lemma we consider the rearrangement of $Tf_x(x)$. This makes it different from the other improved O'Neil lemma in [FM3], which estimated the rearrangement of $Tf(x)$ for a fixed function $f$, not depending on $x$. The proof of Lemma 1, postponed to the Appendix, is based on the proof of Lemma 9 in [FM3], however it is  a bit more streamlined, and it  contains some other minor improvements (see Remark 4 after the proof). Note that Lemma 9 in [FM3] was itself an improvement of Lemma 2 in [FM1], which gave a version of the original O'Neil lemma (see [ON]) for measure spaces. \\

In order to apply the above Lemma \ref{l} to the proof of our main theorem, we also need the following lemma regarding the rearrangement of the sum of two functions whose supports are mutually disjoint.
\bigskip
\begin{lemma}\label{l0}
Let $f_1,f_2:\ M\rightarrow [-\infty,\infty]$ be measurable functions. Suppose that the supports of $f_1$ and $f_2$ are mutually disjoint, $\mu(\textup{supp}f_1)=z$ and \be |f_1|\geq ||f_2||_{\infty} \ \ \ \mu-a.e.\ x\in \{x: f_1(x)\ne 0\}\nt\ee
Then we have
  \be (f_1+f_2)^*(u)=\begin{cases}f_1^*(u) &\text{if} \;\; 0<u< z\\ f_2^*(u-z)& \text{if} \;\;u> z\end{cases},\label{l0a} \ee 
and 
\be (f_1+f_2)^{**}(u)=f_1^{**}(u)\quad \text{for}\ 0<u\le z.\label{I0a1}\ee
\end{lemma}

\ni{\bf Proof of lemma \ref{l0}:} 
Given the assumptions on $f_1,f_2$ we have 
\be m_{f_1+f_2}(s)=m_{f_1}(s)+m_{f_2}(s)\label{I0b1}\ee
and
  \be \begin{cases}m_{f_2}(s)=0 & \text{if}\; m_{f_1}(s)<z\\
m_{f_1}(s)=z & \text{if} \; m_{f_2}(s)>0.\end{cases}\label{l0b} \ee 

For $u<z$, by \eqref{I0b1}, \eqref{l0b} we have $m_{f_1+f_2}(s)\le u$ if $m_{f_1}(s)\le u$. It is also clear that  $m_{f_1}(s)\le u$ whenever $m_{f_1+f_2}(s)\le u$. So $m_{f_1}(s)\leq u$ if and only if $m_{f_1+f_2}(s)\leq u$. We get
\be \begin{aligned}
(f_1+f_2)^*(u)=\inf\{s\geq 0:m_{f_1+f_2}(s)\le u\}=\inf\{s\geq 0: m_{f_1}(s)\le u\}=f_1^*(u). 
\end{aligned}
  \nt\ee
  
 Let $u>z$. If $m_{f_2}(s)\le u-z$ then $m_{f_1+f_2}(s)\le u$. We will show that $m_{f_1+f_2}(s)\le u$ implies $m_{f_2}(s)\le u-z$.  
 Suppose there exists $s\geq 0$ such that $m_{f_1+f_2}(s)\le u$ and $m_{f_2}(s)>0$. Then by \eqref{l0b} we have $m_{f_1}(s)=z$ and hence 
 \be m_{f_2}(s)=m_{f_1+f_2}(s)-m_{f_1}(s)\le u-z.\nt\ee Therefore,
 \be \begin{aligned}
(f_1+f_2)^*(u)=\inf\{s\geq 0: m_{f_1+f_2}(s)\le u\}=\inf\{s\geq 0: m_{f_2}(s)\le u-z\}=f_2^*(u-z).
\end{aligned}
  \nt\ee 
 Lastly, \eqref{I0a1} holds by \eqref{l0a} and the definition of the maximal function.
 \hspace*{\fill} \QEDopen\\


 As a consequence of Lemma \ref{l}, if $f_x=f$ for all $x\in N$, the O'Neil estimate takes the following form:
\be (Tf)^{**}(t)\le C_0t^{-\frac{1}{\beta}} \int_0^{t}{f}^*(u)du+\int_t^{\infty}k_1^*(u){f}^*(u)du.  \nt\ee 
Under the hypothesis that $m(f,s)<\infty$ for $s>0$, we are able to define the O'Neil functional $U$ as follows:
  \be Uf(t)=C_0t^{-\frac{1}{\beta}} \int_0^{t}{f}^*(u)du+\int_t^{\infty}k_1^*(u){f}^*(u)du  \nt\ee 
  where $C_0$ is the constant in the improved O'Neil lemma \ref{l}.
 \par\medskip
  We now  state an Adams inequality due to Fontana and Morpurgo ([FM3, Corollary 2]) in terms of the O'Neil functional. Although the original theorem in their paper is not stated in this form, it is clear from the proof in [FM3] that everything also works for the O'Neil functional instead of the original operator. This result plays a crucial role in the proof of our main result.
  
  \bigskip
 \ni \textup{\textbf{Theorem A ([FM3, Corollary 2])}} \emph{Suppose $\nu(N)<\infty$, $\mu(M)<\infty$, and that
  \be k_1^*(t)\leq A^{\frac{1}{\beta}}t^{-\frac{1}{\beta}}\big(1+H(1+|\log t|)^{-\gamma}\big),\ \ \ 0<t\leq\mu(M)\label{t0a}\ee
  \be k_2^*(t)\leq Bt^{-\frac{1}{\beta}}.\ \ \ 0<t\leq \nu(N)\label{t0b}\ee
  Then there exists a constant $C=C(\b,\gamma,A,B,H)$ such that for each $f\in L^{\beta'}(M)$ with $||f||_{\beta'}\leq 1$, with $\beta^{-1}+(\beta')^{-1}=1$, 
  \be \int_0^{\nu(N)} \exp\bigg[\frac{1}{A}\big(Uf(t)\big)^{\beta}\bigg]dt\leq C\big(\nu(N)+\mu(M)\big).\label{t0c}\ee}

\section{Proof of the inequalities in Theorem \ref{T1}}

\ \ Let us start with following lemma from [FM2], in order to clarify the equivalence between exponential inequalities on sets $\{x\in \R^n\ :\ |Tf(x)|\geq 1\}$ and regularized exponential inequalities over $\R^n$. 

\bigskip
\ni\textup{\bf Lemma A ([FM2, Lemma 9])} \emph{Let $(N,\nu)$ be a measure space and $1<p<\infty$, $a>0$. Then for every $u\in {L}^p(N)$ we have
\be \int_{\{|u|\geq 1\}}e^{a|u|^{p'}}dx-e^a||u||^p_p\leq \int_{N} \bigg(e^{a|u|^{p'}}-\sum_{k=0}^{\lceil p-2\rceil}\frac{a^k|u|^{kp'}}{k!}\bigg)dx\leq \int_{\{|u|\geq 1\}}e^{a|u|^{p'}}dx+e^a||u||^p_p \label{l3a}\ee
and also
\be \int_{\{|u|\geq 1\}}\frac{e^{a|u|^{p'}}}{1+|u|^{p'}}dx-e^a||u||^p_p\leq \int_{N} \frac{e^{a|u|^{p'}}-\sum_{k=0}^{\lceil p-2\rceil}\frac{a^k|u|^{kp'}}{k!}}{1+|u|^{p'}}dx\leq \int_{\{u\geq 1\}}\frac{e^{a|u|^{p'}}}{1+|u|^{p'}}dx+e^a||u||^p_p. \label{l3b}\ee
In particular, the following three inequalities are equivalent:\\
\be \int_{N} \frac{\exp_{\lceil p-2\rceil}{[a|u|^{p'}}]}{1+|u|^{p'}}dx\leq C||u||^p_p,\label{l3cc}\ee
\be \int_{\{|u|\geq 1\}}\frac{e^{a|u|^{p'}}}{1+|u|^{p'}}dx\leq C||u||^p_p,\label{l3dd}\ee
 \be \int_{E} \frac{e^{a|u|^{p'}}}{1+|u|^{p'}}dx\leq C(||u||^p_p+|E|)\label{l3ee}\ee
for every measurable set $E$ with finite measure.}

\bigskip
In order to prove \eqref{1b}, it is enough to show that 
 \be\int_{\{|Tf|\geq 1\}}\frac{\exp\bigg[\dfrac{1}{A_g}|Tf|^{\frac{n}{n-\alpha}}\bigg]}{1+|Tf|^{\frac{ n}{n-\alpha}}}dx \leq C||Tf||_{n/\alpha}^{n/\alpha}. \label{1c}\ee

 
 Let \be t_0=\big|\big\{x\ :\ |Tf|\geq 1\big\}\big|.\nt\ee Note that by this definition we have $(Tf)^*(t)\geq 1$ for $0<t < t_0$ and $(Tf)^*(t)<1$ for $t>t_0$.\par
 Now we will show that \eqref{1c} is equivalent to 
 \be \int_{0}^{t_0}\frac{\exp\bigg[\dfrac{1}{A_g}\big((Tf)^*(t)\big)^{\frac{n}{n-\alpha}}\bigg]}{1+\big((Tf)^*(t)\big)^{\frac{ n}{n-\alpha}}}dt \leq C||Tf||_{n/\alpha}^{n/\alpha}. \label{1d}\ee
 Let us denote the rearrangement of $Tf$ with respect to a measurable set $E$ as $(Tf)_E^*$ and its corresponding maximal function as $(Tf)_E^{**}$. Clearly 
\be (Tf)_E^*(t)=\big((Tf)\chi_E\big)^*(t),\qquad 0<t\le |E|.\label{r16}\ee
 Let 
  \be F(z)=\dfrac{e^{\frac{1}{A_g}z^{n/(n-\a)}}}{1+z^{n/(n-\a)}}\label{r15}\ee
  and $E=\big\{x\ :\ |Tf|\geq 1\big\}$, then the LHS of \eqref{1c} can be written as
 \be \ba \int_{E}F(|Tf|)dx=\int_{0}^{t_0}F((Tf)_E^*(t))dt =\int_{0}^{t_0}F((Tf)^*(t))dt,\ea\label{r17}\ee
The first equality holds since for $F$ non-negative and measurable on $[0,\infty)$, for $g$ measurable on $\R^n$ and if $E$ is a level set of $g$, we have $ \int_{E}F\circ |g|dx= \int_{0}^{|E|}F\circ g_E^*dt$ (see for example [K, Theorem 1.1.1]).
To prove the second equality, note that
\be |(Tf)\chi^{}_{E}(x)|\geq ||(Tf)\chi^{}_{E^c}||_\infty,\qquad\text{for a.e.}\ x\in E,\nt\ee
hence by Lemma \ref{l0} we get 
\be (Tf)_{E}^*(t)=(Tf)^*(t),\qquad\text{for}\ 0<t < t_0.\nt\ee
 
\ \ To estimate $(Tf)^*$, we first define $1$-parameter families of sets $E_\tau,\ F_\tau$ (depending on $f$) as follows.

 For $\tau>0$, let $E_\tau$ be the set such that
\be \bc &|E_{\tau}|=\tau\\
&\{x\ :\ |Tf(x)|>(Tf)^*(\tau)\}\subseteq E_{\tau}\subseteq \{x\ :\ |Tf(x)|\geq(Tf)^*(\tau)\}.
 \ec \label{d2}\ee
 In order to show that such $E_\tau$ exists, we denote $V_1=\{x\ :\ |Tf(x)|>(Tf)^*(\tau)\}$ and $V_2=\{x\ :\ |Tf(x)|\geq(Tf)^*(\tau)\}.$ By definition of rearrangement, we have that $\mu(V_1)\le \tau$ and $\mu(V_2)\geq \tau.$ If $\mu(V_2)=\tau$, we take $E_\tau=V_2$. Otherwise, consider the continuous function $g(r)=\mu(V_1)+\mu(B_r\cap (V_2\setminus V_1))$ for $r\geq 0$, where $B_r=B(0,r)$ is the ball centered at $0$ with radius $r$. It is clear that $g(0)=\mu(V_1)\le \tau$, and $g(r)\to\mu(V_2)$ as $r\to\infty$. Since $\mu(V_2)>\tau$, there exists a $r$ such that $g(r)=\tau$, and $E_\tau=V_1\cup(B_r\cap V_2\setminus V_1)$ is a measurable set that satisfies the condition \eqref{d2}.  
 
 \bigskip
 
Similarly, let $F_{\tau}$  be the set such that  
\be \bc &|F_{\tau}|=\tau\\
&\{x\ :\ |f(x)|>f^*(\tau)\}\subseteq F_{\tau}\subseteq \{x\ :\ |f(x)|\geq f^*(\tau)\}.
 \ec \label{d3}\ee
Let \be f_{\tau}=f\chi^{}_{F_{\tau}},\ \ \ f_\tau'=f\chi^{}_{F^c_{\tau}} \nt\ee
\ni and $r(\tau)=(\tau/|B_1|)^{1/n}$ so that \be |B(0,r(\tau))|=\tau.\label{rtau}\ee 

\begin{rk} If $f$ and $K$  are radially decreasing, then both $E_\tau$ and $F_\tau$ are either open or closed balls of volume $\tau$.\end{rk}

Next, for all $x \in E_\tau$ define 
\be W(\tau,x)=\int_\tau^{2\tau}k_1^*(u)(f_\tau'\chi^{}_{B(x,r(\tau))})^*(u-\tau)du \label{d4} \ee
\be M(\tau,x)=|T(f_\tau'\chi^{}_{B^c(x,r(\tau))})(x)|. \label{d5}\ee
\ni Lastly, for fixed $\tau>0$, take the essential supremum in \eqref{d4} and \eqref{d5}, and let
\be W_\tau=\supess_{x\in E_\tau} W(\tau,x) \label{d7} \ee
\be M_\tau=\supess_{x\in E_\tau} M(\tau,x). \label{d8}\ee
We want to point out that for all $\tau>0$ we have $W_\tau<\infty$ and $M_\tau<\infty$, by the fact that $f$ is compactly supported and $||f||_{n/\a}\leq 1$. 
Also, note that for each $x$ and $\tau$,
\be f=f_\tau+f'_\tau=f_\tau+f_\tau'\chi^{}_{B(x,r(\tau))}+f_\tau'\chi^{}_{B^c(x,r(\tau))}\ \ \ \nt \ee
and \be Tf(x)=Tf_\tau(x)+T(f_\tau'\chi^{}_{B(x,r(\tau))})(x)+T(f_\tau'\chi^{}_{B^c(x,r(\tau))})(x).\label{d6}\ee
From now on we will use the following notation:\be q=\frac{n}{\a},\qquad\qquad q'=\frac{n}{n-\a}.\nt\ee
Recall that the O'Neil functional is defined, with $\b=q'$, as follows
\be Uf(t)=C_0t^{-\frac{1}{q'}} \int_0^{t}{f}^*(u)du+\int_t^{\infty}k_1^*(u){f}^*(u)du.  \nt\ee 
Our first step toward a proof of \eqref{1d} is to establish the following estimate:
\be (Tf)^{*}(t)\leq (Tf)^{**}(t)\leq Uf_\tau(t)+W_\tau+M_\tau\ \ \ \ \ \textup{for}\ 0<t\leq\tau.\label{e1}\ee
Recall the definition of $(Tf)_E^*$ in \eqref{r16}, for any measurable set $E$ . The definition of $E_\tau$ implies that 
\be |(Tf)\chi^{}_{E_\tau}(x)|\geq ||(Tf)\chi^{}_{E^c_\tau}||_\infty,\qquad\text{for a.e.}\ x\in E_\tau,\nt\ee
hence we can apply Lemma \ref{l0} to get 
\be (Tf)_{E_\tau}^{**}(t)=(Tf)^{**}(t),\qquad\text{for}\ 0<t\le\tau.\label{TE1}\ee
Let \be f_{x,\tau}=f_\tau+f_\tau'\chi^{}_{B(x,r(\tau))},\nt\ee
note that by the definition of $M_\tau$ in \eqref{d8} and the decomposition of $Tf$ in \eqref{d6},
\be |(Tf)\chi^{}_{E_\tau}(x)|\le |(Tf_{x,\tau})\chi^{}_{E_\tau}(x)|+M_\tau,\qquad \text{for} \ x\in E_\tau.\nt\ee
Due to subadditivity of $(\cdot)^{**}$ (see [BS, Chapter 2 inequality (3.12)]) and \eqref{TE1}, we have
\be (Tf)^*(t)\le(Tf)^{**}(t)= (Tf)_{E_\tau}^{**}(t)\le (Tf_{x,\tau})_{E_\tau}^{**}(t)+M_\tau,\qquad \ 0<t\le\tau.\label{TE2}\ee
Therefore in order to prove \eqref{e1} it is enough to show the following:
\be  (Tf_{x,\tau})_{E_\tau}^{**}(t)\leq Uf_\tau(t)+W_\tau\qquad \text{for}\  0<t\le \tau.\label{e101}\ee
In other words, we only need to show the rearrangement of $Tf_{x,\tau}(x)$ over the set $E_\tau$ satisfies \eqref{e101}. Let us apply the improved O'Neil Lemma (Lemma \ref{l}) with $N=E_\tau,\ M=\R^n,\ \b=q',\ f_x=f_{x,\tau}$ , and
\be \ol{f}=|f_\tau+f_\tau'|=|f| \label{e2}\ee
so that \eqref{l1} and \eqref{l2} hold. For $x\in E_\tau,$ let 
\be T'f(x)=Tf_{x,\tau}(x)=\int_{\R^n}K(x-y)(f_\tau+f_\tau'\chi_{B(x,r(\tau))})(y)dy. \label{e3}\ee
By Lemma 9 in [FM1], we have that \eqref{A1}, \eqref{A3} implies the condition \eqref{l1}, with $\b=q'$, in Lemma~\ref{l}. From now on we will use $k^*$ to denote $k_1^*$ since $k(x,y)=K(x-y)$ is a convolution kernel. We obtain
\be\ba (T'f)_{E_\tau}^{**}(t)&\leq C_0t^{-\frac{1}{q'}}\int_0^t {f}^*(u)du+\supess_{x\in E_\tau}\int_t^{\infty}k_1^*(x,u)f_{x,\tau}^*(u)du\cr
&=C_0t^{-\frac{1}{q'}}\int_0^t {f}^*(u)du+\supess_{x\in E_\tau}\int_t^{\infty}k^*(u)f_{x,\tau}^*(u)du.\ea\label{e4}\ee
By definition of $F_\tau, f_\tau$ and \eqref{e2}, we apply Lemma \ref{l0} to get
\be {f}^*(u)=(f_\tau+f_\tau')^*(u)=f_\tau^*(u)\ \  \text{if} \;\; 0<u<\tau\label{e5} \ee 
and
\be f_{x,\tau}^*(u)=(f_\tau+f_\tau'\chi_{B(x,r(\tau))})^*(u)=\begin{cases}f_\tau^*(u) &\text{if} \;\; 0<u<\tau\\ \big(f_\tau'\chi_{B(x,r(\tau))}\big)^*(u-\tau)& \text{if} \;\;u>\tau.\end{cases}\label{e6} \ee 
Therefore, \eqref{e4} can rewritten as
\be \ba(T'f)_{E_\tau}^{**}(t)&\leq C_0t^{-\frac{1}{q'}}\int_0^t f_\tau^*(u)du+\int_t^{\tau}k^*(u)f_\tau^*(u)du\cr
&+\supess_{x\in E_\tau}\int_\tau^{2\tau}k^*(u)\big(f_\tau'\chi_{B(x,r(\tau))}\big)^*(u-\tau)du\cr
&=Uf_\tau(t)+W_\tau.\ea \label{e7}\ee
Hence \eqref{e101} is proved and \eqref{e1} follows. 

\bigskip
Next  we consider the following inequality (also in [MS2, the inequality below (4.7)], with slightly different form)
 \be (a+b)^p\leq \lambda^{1-p}a^p+(1-\lambda)^{1-p}b^p\qquad\quad   a,b\geq0,\;0<\lambda<1,\; p>1\label{1*}, \ee  
which can be proved by writing $a+b$ as $(a\lambda^{-1/p'})\lambda^{1/p'}+(b(1-\lambda)^{-1/p'})(1-\lambda)^{1/p'}$ and apply Holder inequality.
Then we use estimation \eqref{e1} and apply the above inequality
to the integrand in \eqref{1d}, with \be \ p=q',\ a=(Uf_\tau)^*(t),\ \textup{and}\ b=W_\tau+M_\tau.\nt\ee We get
\be\ba 
&\frac{\exp\bigg[\dfrac{1}{A_g}\big((Tf)^*(t)\big)^{q'}\bigg]}{1+\big((Tf)^*(t)\big)^{q'}} \leq C\frac{\exp\bigg[\dfrac{1}{A_g}\big(Uf_\tau(t)+W_{\tau}+M_{\tau}\big)^{q'}\bigg]}{1+\big(Uf_\tau(t)+W_{\tau}+M_{\tau}\big)^{q'}}  \cr
&\leq C \frac{\exp\bigg[\dfrac{(1-\lambda)^{1-q'}}{A_g}\big(W_{\tau}+M_{\tau}\big)^{q'}\bigg]}{1+\big(W_{\tau}+M_{\tau}\big)^{q'}}\cdot{\exp\bigg[\dfrac{\lambda^{1-q'}}{A_g}\big(Uf_\tau(t)\big)^{q'}\bigg]}.
\ea
\label{main1}\ee
To get the first inequality in \eqref{main1}, let  $F(z)$ be defined as in \eqref{r15}. Note that $F(z)\geq C>0$ on $[0,\infty)$ and is increasing in $z$ for $z\geq 1+A_g^{(n-a)/n}$. Also recall that for $0<t < t_0$ we have $(Tf)^*(t)\geq 1$. We consider two cases. If $1\le (Tf)^*(t)\le 1+A_g^{(n-a)/n}$, then we have $F((Tf)^*(t))\le C$ and $F(Uf_\tau(t)+W_{\tau}+M_{\tau})\geq C>0$. If  $(Tf)^*(t)\geq 1+A_g^{(n-a)/n}$, then by \eqref{e1} and the fact that $F(z)$ is increasing, the first inequality follows.\par 

 Let $t_1>0$ be the number such that \be \dfrac{\int_0^{t_1}f^*(u)^{q}du}{||f||_{q}^{q}}=\frac{1}{4}\nt\ee and \be\e_\tau=\min\left\{\dfrac{1}{4},\ \dfrac{\int_0^{\tau}f^*(u)^{q}du}{||f||_{q}^{q}}\right\}.\nt\ee
We estimate \eqref{main1} using the following two lemmas. The first one is an integral estimate (essentially the Adams inequality):
\begin{lemma}\label{lemmaI2}
If we define 
\be I_2(\tau,t,\lambda)={\exp\bigg[\dfrac{\lambda^{1-q'}}{A_g}\big(Uf_\tau(t)\big)^{q'}\bigg]}, \qquad \tau>0,\ t>0,\ \lambda>0,\nt\ee
then
\be \int_0^{\tau}I_2(\tau,t,\e_\tau)dt\le C\tau,\qquad 0<\tau\le t_1. \label{c1bbb}\ee
\end{lemma}
\ni {\bf Proof of Lemma \ref{lemmaI2}:} First note that when $\tau\le t_1$, we have 
\be \e_\tau=\dfrac{\int_0^{\tau}f^*(u)^{q}du}{||f||_{q}^{q}}.\nt\ee
If we let  \be \wt{f}:=\frac{f_{\tau}}{\e_{\tau}^{1/q}}, \nt\ee then we have that $\wt{f}$ has measure of support $\mu(\textup{supp} \wt{f})\leq \tau$ with \be ||\wt{f}||_q\leq1,\nt\ee
also assumption \eqref{A1} implies the estimate \eqref{l1} on $k^*$, with $\b=q'$ (See [FM1, Lemma~9]). That is, conditions \eqref{t0a} and \eqref{t0b} are satisfied. Therefore by Theorem~A, the Adams inequality for the O'Neil functional, we obtain \eqref{c1bbb}.\hfill\QEDopen\\
\par
The estimation for $I_1(\tau,\lambda)$ which is stated in the following lemma is essential for the rest of the proof. Let us assume the lemma for now, and its proof will be given in sections 4-5.  
 \begin{lemma}\label{le1}
Let $0<\tau\leq t_0$, $0\le\lambda<1$. Define
\be I_1(\tau,\lambda)=\frac{\exp\bigg[\dfrac{(1-\lambda)^{1-q'}}{A_g}\big(W_{\tau}+M_{\tau}\big)^{q'}\bigg]}{1+\big(W_{\tau}+M_{\tau}\big)^{q'}}.\nt\ee
Then there exists constant $C>0$ such that
\be I_1(\tau,\e_\tau)\leq \frac{C}{\tau}||Tf||^{q}_{q} \label{l1a} \ee
where $C=C(n,\a,K).$
\end{lemma}
Assuming Lemma \ref{le1}, let $ \tau_0=\min\{t_0,t_1\}$. To prove \eqref{1d} it is enough to show that
\be \int_{0}^{\tau_0}\frac{\exp\bigg[\dfrac{1}{A_g}\big((Tf)^*(t)\big)^{q'}\bigg]}{1+\big((Tf)^*(t)\big)^{q'}}dt\leq C||Tf||_q^q \label{c2a} \ee
and then show that if $t_1<t_0$,
\be \int_{t_1}^{t_0}\frac{\exp\bigg[\dfrac{1}{A_g}\big((Tf)^*(t)\big)^{q'}\bigg]}{1+\big((Tf)^*(t)\big)^{q'}}dt\leq C||Tf||_q^q. \label{c2b} \ee
To prove \eqref{c2a}, we take $\tau=\tau_0$ in \eqref{l1a} and \eqref{c1bbb} to get
\be I_1(\tau_0,\e_{\tau_0})\leq \frac{C}{\tau_0}||Tf||^{q}_{q}\qquad\text{and}\qquad \int_0^{\tau_0}I_2(\tau_0,t,\e_{\tau_0})dt \leq C\tau_0. \label{main4} \ee
Therefore, using \eqref{main1} it is immediate that
\be \int_{0}^{\tau_0}\frac{\exp\bigg[\dfrac{1}{A_g}\big((Tf)^*(t)\big)^{q'}\bigg]}{1+\big((Tf)^*(t)\big)^{q'}}dt\leq\int_0^{\tau_0}C I_1(\tau_0,\e_{\tau_0})I_2(\tau_0,t,\e_{\tau_0})dt\le C||Tf||^{q}_{q}\label{main6}\ee
and \eqref{c2a} follows.\\

Next, to show \eqref{c2b}, we take $\tau=t$ for $t_1\le t\le t_0$, and $\lambda=\frac{1}{8}$ in the definition of $I_2$ in Lemma \ref{lemmaI2}. Then by the definition of the O'Neil operator and the fact that the support $f_t$ has measure less than or equal $t$, \be Uf_t(t)=C_0t^{-\frac{1}{q'}}\int_0^t f_t^*(u)du \leq C||f_t||_q\leq C.\nt\ee
So we have \be I_2\bigg(t,t,\frac{1}{8}\bigg)\le C.\label{main7}\ee
Since $t_1\le t\le t_0$ and $\e_{t_1}=\frac{1}{4}$, by definition $\e_t=\frac{1}{4}$. Take $\theta=(6/7)^{\frac{\a}{n-\a}}<1$. Hence
\be \ba I_1\bigg(t,\frac{1}{8}\bigg)&= \frac{\exp\bigg[ \dfrac{(7/8)^{1-q'} }{A_g}\big(W_t+M_t\big)^{q'}\bigg]}{1+(W_t+M_t)^{q'}}=\frac{\exp\bigg[(7/6)^{1-q'} \dfrac{(3/4)^{1-q'} }{A_g}\big(W_t+M_t\big)^{q'}\bigg]}{1+(W_t+M_t)^{q'}}\cr
&=\frac{\left(\exp\bigg[ \dfrac{(3/4)^{1-q'} }{A_g}\big(W_t+M_t\big)^{q'}\bigg]\right)^\theta}{1+(W_t+M_t)^{q'}}\le \Bigg(\frac{\exp\bigg[ \dfrac{(3/4)^{1-q'} }{A_g}\big(W_t+M_t\big)^{q'}\bigg]}{1+(W_t+M_t)^{q'}}\Bigg)^\theta\cr 
&=I_1^{\theta}\bigg(t,\frac{1}{4}\bigg)=I_1^{\theta}(t,\e_t)\le \frac{C}{t^\theta}||Tf||^{\theta q}_q.
\ea
\label{main8}\ee
Using \eqref{main7}, \eqref{main8} we get
\be \ba
\int_{t_1}^{t_0}CI_1\bigg(t,\frac{1}{8}\bigg)I_2\bigg(t,t,\frac{1}{8}\bigg)dt&\leq C\int_{t_1}^{t_0}\frac{1}{t^\theta}||Tf||^{\theta q}_qdt\le Ct_0^{1-\theta}||Tf||^{\theta q}_{q}\cr &\leq C||Tf||^{(1-\theta){q}}_{q}||Tf||^{\theta q}_{q}=C||Tf||_q^q,
\ea
 \label{c2e}\ee
where the last inequality is by the fact that \be ||Tf||_q^q\geq t_0\label{c2f}\ee since $(Tf)^*(t)\geq 1$ for all $t< t_0$, by the definition of $t_0$. So \eqref{c2b} follows from \eqref{main1}.
In order to complete the proof of Theorem \ref{T1}, we are left to prove Lemma \ref{le1}.

\section{Proof of Lemma \ref{le1}}
It is enough to show that 
\be \frac{\exp\bigg[\dfrac{(1-\e_\tau)^{-\frac{\a}{n-\a}}}{A_g}(W(\tau,x_2)+M(\tau,x_1))^{\frac{n}{n-\a}}\bigg]}{1+(W(\tau,x_2)+M(\tau,x_1))^{\frac{n}{n-\a}}}\leq C\frac{||Tf||_q^q}{\tau}\label{le1aaa}\ee
for all $x_1,\ x_2\in E_\tau.$
Now let us state the following key lemma in [MS1]-[MS3], [LTZ].
\begin{lemma}\label{lo}
Let $q>1 $. Given any sequence $\ds a=\{a_k\}_{k\geq 0}$, let \be ||a||_1=\aa |a_k|,\ ||a||_q=\big(\aa |a_k|^q\big)^{1/q}\label{lo1}\ee and define \be \mu_d(h)=\inf \{\aa |a_k|^{q}e^{qk}: \ ||a||_1=h,\ ||a||_q \leq 1\}.\nt\ee Then for $h>1$, we have \be C_1(q)\frac{\exp\big[{qh^{q'}}\big]}{h^{ q'}}\leq \mu_d(h)\leq C_2(q) \frac{\exp\big[{qh^{q'}}\big]}{h^{ q'}}.\label{lo2}\ee\end{lemma}


 As a consequence of  the above optimal growth lemma, we deduce that for any $q>1$ and any $\mu>0,\,h>1$ there is $C=C(q)$ such that 
for any sequence $\{a_k\}$ satisfying 
\be\sum_{k=0}^\infty |a_k|=h \qquad \sum_{k=0}^\infty|a_k|^q\le \mu\label{le1a}\ee we have
\be \frac{\exp\big[{q\mu^{1-q'}h^{q'}}\big]}{h^{q'}}\le C\mu^{-q'}\sum_{k=0}^\infty |a_k|^q e^{qk}.\label{le1b}\ee
The next task is to find a number $h_1$, depending on $f$ and $x_1$, and a sequence $a=\{a_k\}$, also depending on $f$ and $x_1$, such that
  \be q^{-\frac 1{q'}}A_g^{-\frac1 {q'}} (W(\tau,x_2)+M(\tau,x_1))\le h_1\label{le1c}\ee 
\be \sum_{k=0}^{\infty}|a_k|=h_1,\qquad\quad \sum_{k=0}^{\infty}|a_k|^q\le 1-\e_\tau,\qquad\quad \sum_{k=0}^{\infty}|a_k|^qe^{qk}\le \frac {C}{\tau}||Tf||_q^q.\label{le1d} \ee  
Clearly \eqref{le1aaa} follows from \eqref{le1a}-\eqref{le1d}, with $\mu=1-\e_\tau\geq \frac{3}{4}$ and $h=h_1$.\par
From now on, throughout the proof of Lemma \ref{le1}, we fix $0<\tau\leq t_0$, and $x_1,x_2\in E_\tau$ as defined in \eqref{d2}. First let us introduce some notation. Recall that $r(\tau)$ is the number such that $|B(0,r(\tau))|=\tau$. Define for each $j=0,1,2...$
    \be r_j= r(\tau) e^{\frac{q}{n}j},\qquad D_j=B(x_1,\ r_j)  \nt\ee 
  \be \a_j=||f_\tau'\chi_{D_{j+1}\setminus D_{j}}^{}||_q,\qquad \a_{-1}=||f_\tau'\chi_{D_0}||_q\nt\ee \be \ab_j=\max{\{\a_{-1},\a_0,...,\a_j\}},\qquad  \b_j=||f_\tau'\chi_{D^c_{j}}^{}||_q.  \nt\ee 
%
Notice that for any $j$
$$\a_j\le\b_j\le1.$$

 Clearly $\b_j$ is decreasing, and it vanishes for all $j$ large enough, since $f$ has compact support. In particular, there is an integer $N$ so that
   \be \supp f\subseteq D_N=B(x_1,r_N).  \nt\ee 
  Now we are ready to state the main estimates on $M(\tau,x_2)$ and $W(\tau,x_1)$:
\begin{prop}\label{cl}
There exist constants $C_2,C_3$ independent of $f$ and an integer $J$ such that
\be q^{-\frac 1{q'}}A_g^{-\frac1 {q'}} (W(\tau,x_2)+M(\tau,x_1))\le\sum_{j=0}^{J}\a_j+C_2\ab_J+C_2\b_J \label{cl1}\ee
and
\be \sum_{j=0}^{J} \a_j^qe^{qj}+\b_J^qe^{qJ}\leq \frac{C_3}{\tau}||Tf||_q^q.\label{cl2}\ee
\end{prop}
The proof of Proposition \ref{cl} will be given in section 5. Assuming the proposition, we now show how to derive \eqref{l1a}, and hence finish the proof of Lemma~\ref{le1}, using \eqref{le1a}-\eqref{le1d} together with Proposition \ref{cl}. \par
Our goal is to find a number $h_1$ and a sequence $a=\{a_k\}$ that satisfies \eqref{le1c} and \eqref{le1d}. 
Let \be h_1=\sum_{j=0}^J \a_j+C_2\,\ab _J+C_2\b_J.\label{le1e}\ee Clearly, by Proposition \ref{cl}, we have \be q^{-\frac 1{q'}}A_g^{-\frac1 {q'}} (W(\tau,x_2)+M(\tau,x_1))\le h_1.\label{le1f}\ee

Let $J^*$ be the smallest integer such that $\a_{J^*}=\ab_{J}$. It is clear that $J^*\le J$. To construct the sequence $a$ that satisfies \eqref{le1d}, let us first define $N_i,\ i=1,...,4$ as follows:
\be
 \ba
 N_1&=J^*\cr
  N_2&=N_1+\lceil(1+C_2)^{q'}\rceil\cr
  N_3&=N_2+J-1-J^*\cr
  N_4&=N_3+\lceil(1+C_2)^{q'}\rceil.
  \ea
  \label{le1g}\ee
  Let $a=\{a_k\}$ be the following:
  \be
 a_k= \begin{cases}
  \a_{k-1}\;\text{if}\ J^*\ne -1; \ 0 \ \text{if}\ J^*=-1 \; &\text{for}\ \ k=0,...,N_1\cr
 \dfrac{(1+C_2)\ab_{J}}{N_2-N_1}\;\text{if}\ J^*\ne -1; \ \dfrac{C_2\ab_{J}}{N_2-N_1} \ \text{if}\ J^*=-1 \ &\text{for}\ \ k=N_1+1,...,N_2\cr
 \a_{k-N_2+J^*}\;\;&\text{for}\ \ k=N_2+1,...,N_3\cr
  \dfrac{\a_{J}+C_2\b_{J}}{N_4-N_3}\;\; &\text{for}\ \ k=N_3+1,...,N_4. \cr
  \ec
\label{le1h}\ee
  With this definition of $a_k$ we have \be||a||_1=\sum_{k=0}^{N_4}|a_k|=\sum_{k=0}^{N_4}a_k=h_1.\label{le1i}\ee
 If $J^*\ne -1$,
  \be\ba
\sum_{k=0}^{N_4}|a_k|^q&=\sum_{k=0}^{N_1}\a_{k-1}^q+\sum_{k=N_1+1}^{N_2}\left(\frac{(1+C_2)\ab_{J}}{N_2-N_1}\right)^q+\sum_{k=N_2+1}^{N_3}\a_{k-N_2+J^*}^{q}\cr
&+\sum_{k=N_3+1}^{N_4}\left(\frac{\a_{J}+C_2\b_{J}}{N_4-N_3}\right)^q\cr
&\leq \sum_{k=0}^{J^*-1}\a_{k-1}^q+\ab^q_{J}+\sum_{k=J^*+1}^{J-1}\a_k^q+\b^q_{J}\cr
&=||f_\tau'||_q^q\leq (1-\epsilon_\tau)||f||^q_q\leq 1-\epsilon_\tau.
\ea
  \label{le1j}\ee
 Likewise for $J^*=-1$,
 \be\ba
\sum_{k=0}^{N_4}|a_k|^q&=\sum_{k=N_1+1}^{N_2}\left(\frac{C_2\ab_{J}}{N_2-N_1}\right)^q+\sum_{k=N_2+1}^{N_3}\a_{k-N_2+J^*}^{q}\cr
&+\sum_{k=N_3+1}^{N_4}\left(\frac{\a_{J}+C_2\b_{J}}{N_4-N_3}\right)^q\cr
&\leq \a_{-1}^q+\sum_{k=0}^{J-1}\a_{k}^q+\b^q_{J}\cr
&=||f_\tau'||_q^q\leq (1-\epsilon_\tau)||f||^q_q\leq 1-\epsilon_\tau.
\ea
  \label{le1k}\ee
 
And using \eqref{cl2} in Proposition \ref{cl}, we also have, if $J^*\ne-1$
 \be\ba
\sum_{k=0}^{N_4}|a_k|^qe^{qk}&=\sum_{k=0}^{N_1}\a_{k-1}^qe^{qk}+\sum_{k=N_1+1}^{N_2}\left(\frac{(1+C_2)\ab_{J}}{N_2-N_1}\right)^qe^{qk}+\sum_{k=N_2+1}^{N_3}\a_{k-N_2+J^*}^{q}e^{qk}\cr
&+\sum_{k=N_3+1}^{N_4}\left(\frac{\a_{J}+C_2\b_{J}}{N_4-N_3}\right)^qe^{qk}\cr
&\leq \sum_{k=0}^{J^*-1}\a_{k-1}^qe^{qk}+C\ab^q_{J}e^{q(J^*+C_4)}+\sum_{k=J^*+1}^{J-1}\a_k^qe^{q(k+C_4)}+\b^q_{J}e^{q(J+2C_4)}\cr
&\leq Ce^{2C_4}\left(\a_{-1}^q+\sum_{k=0}^{J}\a_{k}^qe^{qk}+\b_J^{q}e^{qJ}\right)\le Ce^{2C_4}\left(C+\sum_{k=0}^{J}\a_{k}^qe^{qk}+\b_J^{q}e^{qJ}\right)\cr
&=Ce^{2C_4}\left(C\frac{\tau}{\tau}+\sum_{k=0}^{J}\a_{k}^qe^{qk}+\b_J^{q}e^{qJ}\right)\le  \frac{C}{\tau}||Tf||_q^q
\ea
  \label{le1l}\ee
  where $C_4$ in the above inequality is $C_4=\lceil(1+C_2)^{q'}\rceil$, and in the last inequality we used the fact that $\tau\le ||Tf||_q^q$ since $(Tf)^*(t)\geq1$ for $0<t\le \tau< t_0$.
 Similarly for $J^*=-1$,we also have 
 \be\ba
\sum_{k=0}^{N_4}|a_k|^qe^{qk}\leq Ce^{2C_4}\left(\a_{-1}^q+\sum_{k=0}^{J}\a_{k}^qe^{qk}+\b_J^{q}e^{qJ}\right)\le \frac{C}{\tau}||Tf||_q^q.
\ea
  \label{le1m}\ee

Finally, \eqref{le1j}-\eqref{le1m} shows that the sequence $a$ satisfies \eqref{le1c} and \eqref{le1d}. Hence \eqref{l1a} follows and the proof is concluded.\hfill\QEDopen

\section{Proof of Proposition \ref{cl}}
  In the following proof we will set for any measurable function $\phi:\R^m\to \R$
   \be S_j \phi=\phi\chi_{D_j^c}^{}=\phi\chi_{\{|y-x_1|\ge r_j\}}^{}.  \nt\ee 
 With this notation we then have
   \be( S_j-S_{j+1})f_\tau'=f_\tau'\chi_{D_{j+1}\setminus D_j}^{}=f_\tau'\chi_{\{r_j\le |y-x_1|<r_{j+1}\}}^{}  \nt\ee
 and
   \be \a_j=\|( S_j-S_{j+1})f_\tau'\|_q,\qquad \b_j=\|S_jf_\tau'\|_q.  \nt\ee 
 
\def\supp{{\rm supp }}\def\N{{\mathbb N}}
%
Also note that 
 \be f_\tau'=f_\tau'\chi^{}_{B(x_1,r(\tau))}+f_\tau'\chi^{}_{B^c(x_1,r(\tau))}=f_\tau'\chi^{}_{D_0}+S_0f_\tau'.\label{de1}\ee
  For the rest of the proof we assume that 
  \be TS_0f_\tau'(x_1) \geq 0.\label{de2}\ee
  If, on the other hand, \be TS_0f_\tau'(x_1) < 0,\label{de22}\ee
  we replace $T$ by $-T$, and the proof is exactly the same.\par
  We first give some preliminary estimates on $W(\tau,x_2)$ and $M(\tau,x_1).$ We have that 
  \be  \ba W(\tau,x_2)&=\int_\tau^{2\tau}k_1^*(u)(f_\tau'\chi^{}_{B(x,r(\tau))})^*(u-\tau)du\cr
  &\leq C\int_\tau^{2\tau}u^{-\frac{1}{q'}}(f_\tau'\chi^{}_{B(x_2,r(\tau))})^*(u-\tau)du\leq C||(f_\tau'\chi^{}_{B(x_2,r(\tau))})^*||_q\cr&= C||f_\tau'\chi^{}_{B(x_2,r(\tau))}||_q.
  \ea \label{la7}\ee
 Since $f$ is supported in $D_N$, we also have that
 \be \supp\ f_\tau'\chi^{}_{B(x_2,r(\tau))} \subseteq D_N=\bigcup_{j=0}^{N-1}(D_{j+1}\setminus D_j)\cup D_0. \nt\ee
  By the definition of $D_j$ it is clear that $B(x_2,r(\tau))$ can only have nonempty intersection with at most two elements in the set \be \big\{D_0,\ D_{j+1}\setminus D_j,\ \textup{for}\ j=0,1,...,N-1\big\},\nt\ee
  therefore we have \be ||f_\tau'\chi^{}_{B(x_2,r(\tau))}||_q\leq \a_{j_1}+\a_{j_2} \label{la8}\ee
  for some $j_1,j_2\in \{-1,0,1,...,N\}$. Then by the definitions of $\a,\ab,\b$, we have for any $J\in\{0,1,...,N\}$ \be \bc \a_{j}\leq \ab_{J}\ \ \ \text{if}\ J\geq j\\ \a_{j}\leq \b_{J}\ \ \ \text{if}\ J\leq j\ec\ \ \ \ j=j_1,j_2\label{la9}\ee
  so that by combining \eqref{la7},\eqref{la8} and \eqref{la9}, we have
   \be W(\tau,x_2)\leq C\ab_J+C\b_J\label{la2}\ee where $C=C(n,\a,K).$
   Next, recall that  \be M(\tau,x)=|T(f_\tau'\chi^{}_{B^c(x,r(\tau))})(x)|. \ee
  By \eqref{de1} and \eqref{de2}, we can write, for any $J\in\{0,1,...,N\}$
  \be \ba M(\tau,x_1)&=|TS_0f_\tau'(x_1)|={T}S_0f_\tau'(x_1) \cr
  &=\sum_{j=0}^J\big({T}S_jf_\tau'(x_1)-{T}S_{j+1}f_\tau'(x_1)\big)+{T}S_{J+1}f_\tau'(x_1)\cr
  &=\sum_{j=0}^J{T}\big(S_jf_\tau'-S_{j+1}f_\tau'\big)(x_1)+{T}S_{J+1}f_\tau'(x_1).
  \ea
  \nt\ee
  For any integer $j$, we have the estimate
   \be\ba &{T}\big(S_jf_\tau'-S_{j+1}f_\tau'\big)(x_1)\le |T(S_jf_\tau'-S_{j+1}f_\tau')(x_1)|\cr&\le \bigg(\mathop\int\limits_{r_j\le|y|<r_{j+1}}|K(y)|^{q'}dy\bigg)^{1/q'}\|S_jf_\tau'-S_{j+1}f_\tau'\|_q.\ea\label{la3} \ee 
   Using \eqref{A1}, \eqref{A4} and  the inequality $(a+b)^\beta\le a^\b+\b 2^{\b-1}(a^{\b-1}b+b^\b)$ for $\b>1$ (see Adams [A1, inequality (17)] or use mean value theorem) we  get
   \be |K(y)|^{q'}\le |g(y^*)|^{q'}|y|^{-n}+C\min\{|y|^{-n+\delta_1}, |y|^{-n-\delta_2}\}  \label{la4}\ee 
   for some $C>0,\ C=C(n,\alpha,H_1,H_2,B,\delta_1,\delta_2)$.
   Since $r_{j+1}=e^{\frac q n}r_j$, 
    \be \ba&{T}\big(S_jf_\tau'-S_{j+1}f_\tau'\big)(x_1)\le \big(qA_g+C\min\{r_j^{\delta_1}, r_j^{-\delta_2}\}\big)^{1/q'}\a_j\cr
    &\le \left(q^{\frac1{q'}}A_g^{\frac1{q'}}+C\min\{r_j^{{\delta_1/q'}}, r_j^{-{\delta_2/q'}}\}\right)\a_j.\ea\label{la5} \ee 
    Using \eqref{la5}, we then get that
    \be \ba M(\tau,x_1)&\leq \sum_{j=0}^J\left(q^{\frac1{q'}}A_g^{\frac1{q'}}+C\min\{r_j^{{\delta_1/q'}}, r_j^{-{\delta_2/q'}}\}\right)\a_j+{T}S_{J+1}f_\tau'(x_1)\cr    
    &=q^{\frac1{q'}}A_g^{\frac1{q'}} \sum_ {j=0}^{J}\a_j+ C\ab_J \sum_{j=0}^\infty \min\{r_j^{{\delta_1/q'}}, r_j^{-{\delta_2/q'}}\}+ TS_{J+1}f_\tau'(x_1)\cr
    &\le q^{\frac1{q'}}A_g^{\frac1{q'}} \sum_ {j=0}^{J}\a_j+ C\ab_J \sum_{j=0}^\infty (e^{-\frac{q}{n}\frac{\delta_1}{q'}j}+e^{-\frac{q}{n}\frac{\delta_2}{q'}j})+ TS_{J+1}f_\tau'(x_1)\cr
 &= \ q^{\frac1{q'}}A_g^{\frac1{q'}} \sum_ {j=0}^{J}\a_j+ C\ab_J + TS_{J+1}f_\tau'(x_1).
  \ea
  \label{la1}\ee
  \par
  Note that \eqref{la9} and \eqref{la1} are true for any $J\in\{0,1,...,N\}$. The main task now is prove that there exists $J\in\{0,1,...,N\}$ such that
\be TS_{J+1}f_\tau'(x_1)\le  C\ab_J +C\b_J \label{clcl2}\ee
and that \eqref{cl2} holds.
This will be effected by a double stopping time argument, which will simultaneously yield \eqref{cl2} in Proposition \ref{cl}.\par
 Recall that $N$ is an integer such that $\supp f\subseteq D_N$.  Let $J_1\in\{1,...,N\} $ be such that
\begin{numcases}{}
 \b_{j+1}^q\le \left(\avgdj |Tf(x)|dx\right)^q  &  {\text {for  }}  $j=0,...,J_1-1$\label {j1a}\\
 \b_{J_1+1}^q>\left(\avint_{D_{J_1+1}\setminus D_{J_1}} | Tf(x)|dx\right)^q.   \label{j1b}
 \end{numcases}
 If condition \eqref{j1a} is never satisfied we let $J_1=0$, and if \eqref{j1b} is never satisfied let $J_1=N+1$. Next, let $J_2\in\{1,...,N\} $be such that
 \begin{numcases}{}
  {T}S_{j+1}f_\tau'(x_1)\geq \left(\frac{e^{q-1}+1}{2e^{q-1}} \right) {T}S_{j}f_\tau'(x_1) &  {\text {for  }}  $j=0,...,J_2-1$\label {j2a}\\
{T}S_{J_2+1}f_\tau'(x_1)< \left(\frac{e^{q-1}+1}{2e^{q-1}}  \right) {T}S_{J_2}f_\tau'(x_1)  . \label{j2b}
 \end{numcases}
 As in the definition of $J_1$, we let $J_2=0$ if condition \eqref{j2a} is never satisfied , and let $J_2=N+1$  if \eqref{j2b} is never satisfied. \par
 We will first prove \eqref{clcl2}, and hence \eqref{cl1}, in three cases depending on $J_1,J_2$, then we will show that \eqref{cl2} holds with the chosen $J$ in each case.
 
 \bigskip
 \ni{\bf \underline{Case 1:}} $J_2\leq J_1\leq N+1$ and $J_2\ne N+1$.
 
  \bigskip
 \ni{\bf \underline{Case 2:}} $J_2\geq J_1+1$.
 
  \bigskip
 \ni{\bf \underline{Case 3:}} $J_1=J_2=N+1$.
 
  \bigskip
 \ni {\underline{\emph{Proof of \eqref{clcl2} in the case $J_2\leq J_1\leq N+1$ and $J_2\ne N+1$:}}}
 
  \bigskip
 In this case, by \eqref{j2b} we have 
 \be \ba {T}S_{J_2+1}f_\tau'(x_1)&< \left(\frac{e^{q-1}+1}{2e^{q-1}}  \right){T}S_{J_2}f_\tau'(x_1) \cr
 &= \left(\frac{e^{q-1}+1}{2e^{q-1}}  \right) \big({T}S_{J_2+1}f_\tau'(x_1)+{T}(S_{J_2}-S_{J_2+1})f_\tau'(x_1)\big)\cr
 &\leq \left(\frac{e^{q-1}+1}{2e^{q-1}}  \right) \big({T}S_{J_2+1}f_\tau'(x_1)+\big|{T}(S_{J_2}-S_{J_2+1})f_\tau'(x_1)\big|\big)\cr
 &\le \left(\frac{e^{q-1}+1}{2e^{q-1}}  \right) \big({T}S_{J_2+1}f_\tau'(x_1)+C\a_{J_2}\big)\cr
 \ea\label{cl3}\ee
 where the last inequality is by \eqref{la5}. So we have 
 \be {T}S_{J_2+1}f_\tau'(x_1)<C\left(\frac{2e^{q-1}}{e^{q-1}-1}  \right)\a_{J_2}=C\a_{J_2}.\label{cl4}\ee
 Hence by taking $J=J_2$ in \eqref{la1}, we obtain
 \be \ba q^{-\frac 1{q'}}A_g^{-\frac1 {q'}} (W(\tau,x_2)+M(\tau,x_1))&\le\sum_{j=0}^{J_2}\a_j+C\ab_{J_2}+C\b_{J_2}+{T}S_{J_2+1}f_\tau'(x_1)\cr
 &\le\sum_{j=0}^{J_2}\a_j+C\ab_{J_2}+C\b_{J_2}+C\a_{J_2}\cr
 &\le \sum_{j=0}^{J_2}\a_j+C\ab_{J_2}+C\b_{J_2}. \ea\label{cl5}\ee
 Therefore we get \eqref{cl1} with $J=J_2$ and $C_2=C$ in the above inequality.\\
 
  \ni {\underline{\emph{Proof of \eqref{clcl2} in the case $J_2\geq J_1+1$:}}}
  
   \bigskip
 We will need the following lemma to handle this case. Let us state it here, and its proof will be postponed to the Appendix.
 \begin{lemma}\label{lb}
  There is a constant $C_1=C_1(n,\a,K)$ such that for any $J\le N-1$ 
  \be \avgdJ |Tf_\tau(x)|dx\leq C_1\left(\frac{1}{e^{q-1}}\right)^J,\label{lb1}\ee
  
\be \bigg| \avgdJ {T}S_{J+2}f_\tau'(x)dx- {T}S_{J+1}f_\tau'(x_1)\bigg|\le C_1\b_{J+1},\label{lc1}\ee

 \be \bigg|\avgdJ {T}(S_0-S_{J+2})f_\tau'(x)dx\bigg|
\le C_1\ab_{J+1},\label{ld1}\ee
and
\be \bigg|\avgdJ T(f_\tau'\chi_{D_0})(x)dx\bigg|
\le C_1\a_{-1}.\label{ld2}\ee
\end{lemma} 
 %
  Assuming Lemma \ref{lb}, let us first make a reduction. Recall that $0<\tau\le t_0.$ We will assume that \be M(\tau,x_1)\geq\max\{4C_1,1\}. \label{rdt1}\ee where $C_1$ is the constant which is defined in Lemma \ref{lb}. If the above is not true, then we have that $M(\tau,x_1)\le C$, and on the other hand, by \eqref{la7}
\be \ba W(\tau,x_2)\le C||f_\tau'\chi^{}_{B(x_2,r(\tau))}||_q\le C.
  \ea
  \label{la777}\ee
Therefore, $W(\tau,x_2)+M(\tau,x_1)\le C$, and hence 
\be \frac{\exp\bigg[\dfrac{(1-\e_\tau)^{-\frac{\a}{n-\a}}}{A_g}(W(\tau,x_2)+M(\tau,x_1))^{\frac{n}{n-\a}}\bigg]}{1+(W(\tau,x_2)+M(\tau,x_1))^{\frac{n}{n-\a}}}\leq C= C\frac{\tau}{\tau}\le C\frac{||Tf||_q^q}{\tau}\label{rdt2}\ee
which is \eqref{le1aaa}, and the last inequality is by \eqref{c2f}.\\ 
%
 By \eqref{lc1} in Lemma \ref{lb} and recalling that 
 \be f=f_\tau+f_\tau'=f_\tau+f_\tau'\chi^{}_{D_0}+f_\tau'\chi^{}_{D_{J_1+2}\setminus D_0}+f_\tau'\chi^{}_{D_{J_1+2}^c},\nt\ee
  we have
 \be \ba & {T}S_{J_1+1}f_\tau'(x_1)\leq \avgdJo {T}S_{J_1+2}f_\tau'(x)dx+C\b_{J_1+1}\cr
 &\leq \bigg|\avgdJo {T}S_{J_1+2}f_\tau'(x)dx\bigg|+C\b_{J_1+1}\cr
 &=\bigg|\avgdJo\big(Tf-Tf_\tau-T(f_\tau'\chi^{}_{D_0}) -T(S_0-S_{J_1+2})f_\tau'\big)(x)dx\bigg|+C\b_{J_1+1}\cr
 &\le\bigg|\avgdJo Tf(x)dx\bigg|+\bigg|\avgdJo Tf_\tau(x)dx\bigg|+\bigg|\avgdJo T(f_\tau '\chi^{}_{D_0})(x)dx\bigg|\cr&+\bigg|\avgdJo T(S_0-S_{J_1+2})f_\tau'(x)dx\bigg|+C\b_{J_1+1}\cr
 &\le \avgdJo |Tf(x)|dx+\avgdJo |Tf_\tau(x)|dx+C\a_{-1}+C\ab_{J_1+1}+C\b_{J_1+1}
 \ea\label{cl8}\ee
 where the last inequality is by Lemma \ref{lb} \eqref{ld1},\eqref{ld2}.
 To estimate the second integral, note first that by reduction \eqref{rdt1} we have \be M(\tau,x_1)={T}S_0f_\tau'(x_1)\geq 4C_1. \label{cl7}\ee
 Using \eqref{lb1} in Lemma \ref{lb}, and condition \eqref{j2a}, we get
 \be \ba\avgdJo |Tf_\tau(x)|dx&\leq C_1\left(\frac{1}{e^{q-1}}\right)^{J_1} \leq \frac{1}{4}\left(\frac{1}{e^{q-1}}\right)^{J_1}{T}S_0f_\tau'(x_1)\cr&\le \frac{1}{4}\left(\frac{1}{e^{q-1}}\right)^{J_1}\left(\frac{2e^{q-1}}{e^{q-1}+1}\right)^{J_1+1}{T}S_{J_1+1}f_\tau'(x_1)\cr
 &= \frac{1}{4}\left(\frac{2}{e^{q-1}+1}\right)^{J_1}\left(\frac{2e^{q-1}}{e^{q-1}+1}\right){T}S_{J_1+1}f_\tau'(x_1)\le \frac{1}{2}{T}S_{J_1+1}f_\tau'(x_1).
 \ea\label{cl7}\ee
 Hence we have 
 \be \ba  {T}S_{J_1+1}f_\tau'(x_1)\leq  \avgdJo |Tf(x)|dx+\frac{1}{2}{T}S_{J_1+1}f_\tau'(x_1)+C\ab_{J_1+1}+C\b_{J_1+1}.
 \ea\nt\ee
 So the above inequality along with the condition \eqref{j1b} give us 
 \be \ba {T}S_{J_1+1}f_\tau'(x_1)&\leq 2\avgdJo |Tf(x)|dx+2C\ab_{J_1+1}+2C\b_{J_1+1}\cr
 &\le 2C\ab_{J_1+1}+(2C+2)\b_{J_1+1}.\ea\label{cl9}\ee
 By taking $J=J_1$ in \eqref{la1}, we get 
 \be \ba q^{-\frac 1{q'}}A_g^{-\frac1 {q'}} (W(\tau,x_2)+M(\tau,x_1))&\le\sum_{j=0}^{J_1}\a_j+C\ab_{J_1}+C\b_{J_1}+{T}S_{J_1+1}f_\tau'(x_1)\cr
 &\le\sum_{j=0}^{J_1}\a_j+C\ab_{J_1+1}+C\b_{J_1+1}\cr&\le \sum_{j=0}^{J_1}\a_j+C\ab_{J_1}+C\b_{J_1} \ea\label{cl10}\ee
 where the last inequality is by the fact that $\ab_{J+1}\le\ab_J+\a_{J+1}\le  \ab_J+\b_{J+1}\le \ab_J+\b_J$. Therefore we get \eqref{cl1} with $J=J_1$.
 
  \bigskip
 \ni {\underline{\emph{Proof of \eqref{clcl2} in the case $J_1=J_2=N+1$:}}}
 
  \bigskip
In this case, we will simply write the entire series, that is, we will take $J=N$. Since ${T}S_{N+1}f_\tau'(x_1)=0$ we have
\be\ba q^{-\frac 1{q'}}A_g^{-\frac1 {q'}} (W(\tau,x_2)+M(\tau,x_1))&\le\sum_{j=0}^{N}\a_j+C\ab_{N}+C\b_{N}+{T}S_{N+1}f_\tau'(x_1)\cr
 &\le \sum_{j=0}^{N}\a_j+C\ab_{N}.\ea\label{cl12}\ee
 %
  To check \eqref{cl2}, note that we take $J=J_2$ in case 1, $J=J_1$ in case 2 and $J=N$ in case 3. Assume first that $J_1\neq 0$. Then in all the cases we have that \eqref{j1a} is true for all $j\leq J-1$,
 so 
 \be  \ba \sum_{j=0}^{J} \a_j^qe^{qj}+\b_{J}^qe^{qJ}&\leq 3e^{2q}\sum_{j=0}^{J-1} \b_{j+1}^qe^{qj}+3e^{2q}\b_0^q\le C\sum_{j=0}^{J-1} \left(\avgdj |Tf(x)|dx\right)^qe^{qj}+C\cr
 &\le C\sum_{j=0}^{J-1}  \left(\frac{1}{r_j^n} \dj |Tf(x)|dx\right)^qe^{qj}+C\cr
 &\le C\sum_{j=0}^{J-1} r_j^{-qn}\left( \dj |Tf(x)|^qdx\right)\left( \dj dx\right)^{q/q'}e^{qj}+C\cr
 &=C\sum_{j=0}^{J-1} r_j^{-qn}r_j^{(n-\a)q}\left( \dj |Tf(x)|^qdx\right)e^{qj}+C\cr
 &\le C\sum_{j=0}^{J-1}r_j^{-n}e^{qj}\dj |Tf(x)|^qdx+C\le \frac{C}{\tau}||Tf||_q^q
 \ea\label{cl6}\ee
 where in the last inequality we used \eqref{c2f} and also the fact that $r_j^n=r_0^ne^{qj}$ and $\tau=|B(x_1,r_0)|$, so $\tau=Cr_0^n$.\\
 If $J_1=0$, then we just need to check \eqref{cl2} for $J=0$:
 \be  \a_0^q+\b_{0}^q\leq 2\b_0^q\le C=C\frac{\tau}{\tau}\le \frac{C}{\tau}||Tf||_q^q, \label{cl666}\ee
 where the last inequality is by \eqref{c2f}. 
 Proposition \ref{cl} is proved. \hfill\QEDopen \\

\section{ Proofs of the sharpness statements in Theorem~\ref{T1}.}\label{sharpness}
\par 
\ \ \  We will make use of the  extremal family of functions constructed in  [FM2, Section 6], that the authors used to prove the sharpness of the exponential constants in  \eqref{4c} and \eqref{3c}.  Specifically, under the hypothesis that $K$ is $n$-regular, the authors produced a family of compactly supported functions $\psi_{\e,r}\in L^q(B(0,r))$ such that 

\be \max\{||\psi_{\e,r}||_q^{q}\ ,\ ||T\psi_{\e,r}||_{ q}^{ q}\}\le 1\nt\ee

\be\ba |T{\psi}_{\e,r}(x)|^{q'}&\geq A_g\log{\frac{1}{(\e r)^n}}+b_r\bigg(1-\frac{C}{\log{\frac{1}{\e^n}}}\bigg)-C,\qquad |x|\le \e r/2,\ea\label{shp15}\ee

\be ||T\psi_{\e,r}||_{ q}^{q}\le Cr^n(\log{\frac{1}{\e^n}})^{-1},\label{normtf2}\ee

where \be b_r:=\int_{1\le |y|\le r }|K(y)|^{q'}dy\nt\ee
and \be 1\le r^{n}\le \frac{A_g}{2C_4}\bigg(\log{\frac{1}{\e^n}}\bigg).\label{shp13}\ee
Note that by the assumptions \eqref{A1}, \eqref{A4}, we have
\be b_r\le A_g\log r^n+C.\label{sp9}\ee

Note also that for $\e$ small  \eqref{shp15} and \eqref{shp13} imply
\be  |T\psi_{\e,r}(x)|\ge 1,\qquad \forall x\in B_{\e r/2}. \label{TF}\ee

To prove that the exponential constant sharp, i.e. it cannot be replaced by a larger constant,  pick 
$$r^{n}= \frac{A_g}{2C_4}\bigg(\log{\frac{1}{\e^n}}\bigg)$$
and for any fixed $\theta>1$ estimate

\be \ba &\int_{\{|T\psi_{\e,r}|\ge1\}}\frac{\exp\bigg[\dfrac{\theta}{A_g}|T\psi_{\e,r}(x)|^{\frac{n}{n-\alpha}}\bigg]}{1+|T\psi_{\e,r}(x)|^{\frac{ n}{n-\a}}}dx\ge  \int_{B_{\e r/2}}\frac{\exp\bigg[\dfrac{\theta}{A_g}|T\psi_{\e,r}(x)|^{\frac{n}{n-\alpha}}\bigg]}{1+|T\psi_{\e,r}(x)|^{\frac{ n}{n-\a}}}dx\cr&\qquad \geq |B_{\e r/2}|\frac{\exp\bigg[\theta\log{\dfrac{1}{(\e r)^n}}+\dfrac{\theta b_r}{A_g}\bigg(1-\dfrac{C}{\log{\frac{1}{\e^n}}}\bigg)-\theta C\bigg]}{1+A_g\log{\dfrac{1}{(\e r)^n}}+Cb_r-C}\geq \dfrac{C(\e r)^{-(\theta-1)n}}{\log{\dfrac{1}{(\e r)^n}}}\rightarrow\infty\ea\label{S17}\ee
as $\e\to 0^+$, and 
where the last inequality is by the estimate of $b_r$ in \eqref{sp9}.\\

Using exponential regularization, Lemma A, we get, for any $\theta>1$

\be\lim_{\e\to0^+}\int_{\R^n} \frac{\exp_{\lceil\frac n\a -2\rceil}\bigg[\dfrac{\theta}{A_g}|T\psi_{\e,r}(x)|^{\frac{n}{n-\alpha}}\bigg]}{1+|T\psi_{\e,r}(x)|^{\frac{ n}{n-\a}}}dx=+\infty,\ee

which proves the sharpness of the exponential constant.

\medskip
To show the sharpness of the power of the denominator 
take  $r=1$, so that  $b_r=0$. For any fixed $\theta<1$ we have
\be \ba \int_{B_{\e/2}}\frac{\exp\bigg[\dfrac{1}{A_g}|T\psi_{\e,1}(x)|^{\frac{n}{n-\alpha}}\bigg]}{1+|T\psi_{\e,1}(x)|^{\frac{\theta n}{n-\a}}}dx&\geq C\e^n\frac{\exp\bigg[\log{\dfrac{1}{\e^n}}-C\bigg]}{1+\Big(A_g\log{\dfrac{1}{\e^n}}-C\Big)^{\theta}}\cr
&\geq \frac{C}{1+\Big(\log{\dfrac{1}{\e^n}}\Big)^{\theta}}\geq C\Big(\log{\dfrac{1}{\e^n}}\Big)^{-\theta}.\ea\label{S18}\ee
Therefore by the estimation \eqref{normtf2} on the $ q$-th norm of $T\psi_{\e,1}$ we have, for any $\theta<1$, 

$$\lim_{\e\to0^+}\|T\psi_{\e,1}\|_q^{-q}\int_{\R^n} \frac{\exp_{\lceil\frac n\a -2\rceil}\bigg[\dfrac{1}{A_g}|T\psi_{\e,1}(x)|^{\frac{n}{n-\alpha}}\bigg]}{1+|T\psi_{\e,1}(x)|^{\frac{\theta n}{n-\a}}}dx=+\infty.$$

\bigskip
\begin{rk}\label{vectorsharp}{For the vector case the proof of sharpness is almost the same. We take the family of functions as in [FM2, Section 6] and the rest of the proof still
works.}\end{rk}

 \bigskip
\begin{rk}\label{example}{An example where the inequality in Theorem \ref{T1} fails but \eqref{4c} holds.}\end{rk}
\ \ \ For an example that Theorem \ref{T1} cannot hold merely under the assumption that $K$ is a Riesz-like kernel, we can take $0<\a<\frac{n}{2}$ and let $K\in C^1(\R^n\setminus 0)$ be such that
\be K(x)=\bc |x|^{\a-n}\ \ &\text{if}\ |x|\le 1\\ 2|x|^{\a-n}\ \ &\text{if}\ |x|\geq 2.\ec\nt\ee
Note that we have \be 2^{q'}|B_1|\log r^n-C\le b_r\le 2^{q'}|B_1|\log r^n+C.\nt\ee
Choose
\be r^{n}=\frac{A_g}{2C_4}\bigg(\log{\frac{1}{\e^n}}\bigg),\label{shp130}\ee
which satisfies \eqref{shp13}. Hence we have
\be \ba \int_{B_{\e r/2}}\frac{\exp\bigg[\dfrac{1}{|B_1|}|T\psi_{\e,r}(x)|^{\frac{n}{n-\alpha}}\bigg]}{1+|T\psi_{\e,r}(x)|^{\frac{ n}{n-\a}}}dx&\geq |B_{\e r/2}|\frac{\exp\bigg[\log{\dfrac{1}{(\e r)^n}}+\dfrac{ b_r}{|B_1|}\bigg(1-\dfrac{C}{\log{\frac{1}{\e^n}}}\bigg)-C\bigg]}{1+C\log{\dfrac{1}{(\e r)^n}}+b_r-C}\cr
&\geq C\frac{r^{2^{q'}n}}{1+Cr^n}\to\infty\ea\label{S20}\ee
 as $\e\rightarrow 0^+$.\par
 On the other hand, since $K$ is a Riesz-like kernel, the inequality \eqref{4c} [FM2, Theorem 5] holds under the Ruf condition. 

\section{Proof of Corollary \ref{T3}}
Assume that
\be||f||_{n/\alpha}\leq 1.\label{T3a0}\ee
Let $q=n/\a.$ It is enough to show that
\be  \int_{\{|Tf|\geq 1\}}\exp\bigg[\frac{\theta}{A_g}|Tf(x)|^{\frac{n}{n-\alpha}}\bigg]dx\leq \frac{C}{1-\theta}||Tf||_{q}^{q}   \label{T3a}\ee
since \eqref{3c} is then a direct consequence of the exponential regularization Lemma A. 

To show \eqref{T3a}, write
\be \ba
&\exp\bigg[\frac{\theta}{A_g}|Tf(x)|^{\frac{n}{n-\alpha}}\bigg]\cr
&=\frac{\exp\big[\frac{1}{A_g}|Tf(x)|^{\frac{n}{n-\alpha}}\big]}{1+|Tf|^{\frac{n}{n-\a}}}\frac{1+|Tf|^{\frac{n}{n-\a}}}{\exp\big[\frac{1-\theta}{A_g}|Tf(x)|^{\frac{n}{n-\alpha}}\big]}.
\ea
\label{T3c}\ee
Observe that \be \frac{1+y}{e^{(1-\theta)y/{A_g}}} \leq \frac{C}{1-\theta},\ \ \ \text{for} \ y\geq 0\nt\ee
So by Theorem \ref{T1},
\be\ba \int_{\{|Tf|\geq 1\}}\exp\bigg[\frac{\theta}{A_g}|Tf(x)|^{\frac{n}{n-\alpha}}\bigg]dx
&\le  \frac{C}{1-\theta} \int_{\{|Tf|\geq 1\}}\frac{\exp\bigg[\frac{1}{A_g}|Tf(x)|^{\frac{n}{n-\alpha}}\bigg]}{1+|Tf|^{\frac{n}{n-\a}}}dx\cr
&\le \frac{C}{1-\theta}||Tf||_q^q.
\ea
\label{T3d}\ee
Obviously \eqref{T3d} also follows under the more restrictive condition
\be||f||^{pn/\alpha}_{n/\alpha}+||Tf||^{pn/\alpha}_{n/\alpha}\leq 1,\ \ \  p<\infty.\nt\ee
The proof of sharpness is the same as in [FM2]. We use the family of functions $\psi_{\e,r}$  in section \ref{sharpness}, and choose \be r^{n}=\frac{A_g}{2C_4}\bigg(\log{\frac{1}{\e^n}}\bigg).\nt\ee
\section{Proof that Theorem \ref{T1} implies \eqref{4c}}\label{rufinq}
It is enough to show that 
\be \int_{\{|Tf|\geq 1\}}\exp\bigg[\frac{1}{A_g}|Tf(x)|^{\frac{n}{n-\alpha}}\bigg]dx\leq C    \label{T4a}\ee
under the Ruf condition
\be ||f||^{n/\alpha}_{n/\alpha}+||Tf||^{n/\alpha}_{n/\alpha}\leq 1.\nt\ee
Let $\tau=||Tf||_q^q$. Clearly we can assume that $\tau\in(0,1)$. We consider two cases:

 \bigskip
\underline{\bf Case 1:} $\tau\geq 1-({2}/{3})^{q-1}.$

 \bigskip
\underline{\bf Case 2:} $\tau<1-({2}/{3})^{q-1}.$

 \bigskip
\underline{\emph{Proof of \eqref{T4a} in case 1:}} In this case, \be ||f||_q^q\leq 1-\tau\le \big(\frac{2}{3}\big)^{q-1},\label{T4b}\ee
so letting $\wt{f}=f/(\frac{2}{3})^{\frac{q-1}{q}}=\big(\frac{3}{2}\big)^{\frac{q-1}{q}}f$ gives $||\wt{f}||_q^q\le 1$.  We can write
\be \int_{\{|Tf|\geq 1\}}\exp\bigg[\frac{1}{A_g}|Tf(x)|^{\frac{n}{n-\alpha}}\bigg]dx=\int_{\{|Tf|\geq 1\}}\exp\bigg[\dfrac{2}{3A_g}|T\wt{f}(x)|^{\frac{n}{n-\alpha}}\bigg]dx. \label{T4c}\ee
 So by taking $\theta=\dfrac{2}{3}$ in Adachi-Tanaka result, we have 
\be \int_{\{|Tf|\geq 1\}}\exp\bigg[\frac{2}{3A_g}|T\wt{f}(x)|^{\frac{n}{n-\alpha}}\bigg]dx\le \frac{C}{1-2/3}||T\wt{f}||_q^q=3{(\frac{3}{2})^{q-1}C}||T{f}||_q^q\le C.\label{T4d}\ee
Combining \eqref{T4c} and \eqref{T4d} finishes the proof in case 1.

 \bigskip
\underline{\emph{Proof of \eqref{T4a} in case 2:}} In this case, 
\be ||f||_q^q\le 1-\tau \in\bigg(\big(\frac{2}{3}\big)^{q-1},1\bigg).\label{T4e}\ee
Let $p>1$ be such that
\be p(1-\tau)^{\frac{\a}{n-\a}}=1. \nt\ee
Rewrite \eqref{T4a} and apply Holder's inequality,
\be \ba&\int_{\{|Tf|\geq 1\}}\exp\bigg[\frac{1}{A_g}|Tf(x)|^{\frac{n}{n-\alpha}}\bigg]dx= \int_{\{|Tf|\geq 1\}} \frac{\exp\bigg[\dfrac{1}{A_g}|Tf(x)|^{\frac{n}{n-\alpha}}\bigg]}{1+|Tf|^{\frac{n}{(n-\a)p}}}\bigg(1+|Tf|^{\frac{n}{(n-\a)p}}\bigg)dx\cr
&\le \left(\int_{\{|Tf|\geq 1\}} \left(\frac{\exp\bigg[\dfrac{1}{A_g}|Tf(x)|^{\frac{n}{n-\alpha}}\bigg]}{1+|Tf|^{\frac{n}{(n-\a)p}}}\right)^p dx\right)^{\frac{1}{p}} \left(\int_{\{|Tf|\geq 1\}} \bigg(1+|Tf|^{\frac{n}{(n-\a)p}}\bigg)^{\frac{p}{p-1}}dx\right)^{\frac{p-1}{p}}\cr
&=I'\cdot I''.
\ea\label{T4f}\ee
Let \be \wt{f}=\frac{f}{(1-\tau)^{1/q}}\nt\ee
so that $||\wt{f}||_q^q\le 1$. Applying Theorem \ref{T1} gives
\be \ba I'&\le \frac{C}{(1-\tau)^{1/q}}\left(\int_{\{|Tf|\geq 1\}} \left(\frac{\exp\bigg[\dfrac{(1-\tau)^{\frac{\a}{n-\a}}}{A_g}|T\wt{f}(x)|^{\frac{n}{n-\alpha}}\bigg]}{1+|T\wt{f}|^{\frac{n}{(n-\a)p}}}\right)^p dx\right)^{\frac{1}{p}}\cr
&\le C\left(\int_{\{|Tf|\geq 1\}} \frac{\exp\bigg[\dfrac{1}{A_g}|T\wt{f}(x)|^{\frac{n}{n-\alpha}}\bigg]}{1+|T\wt{f}|^{\frac{n}{(n-\a)}}} dx\right)^{\frac{1}{p}}\cr
&\le C||T\wt{f}||_q^{\frac{q}{p}}\le \frac{C}{1-\tau}||Tf||_q^{\frac{q}{p}}\le C||Tf||_q^{\frac{q}{p}}.
\ea\label{T4g}\ee
To estimate $I''$ we start with the following Adachi-Tanaka inequality:
\be \int_{\{|Tf|\geq 1\}}\exp\bigg[\frac{1}{2A_g}|T{f}(x)|^{\frac{n}{n-\alpha}}\bigg]dx\le C||T{f}||_q^q.\label{T4h}\ee

Let \be P=\frac 1{p-1},\ P_0=\bigg\lceil\frac{1}{p-1}\bigg\rceil-1,\ P_1=\bigg\lceil\frac{1}{p-1}\bigg\rceil\nt\ee and define \be F(x):=\bc Tf(x) \ \ &\text{if}\ |Tf(x)|\geq 1\\ 0\ \ &\text{otherwise}.\ec\ee
By the power series expansion of the exponential function, we have that the inequality \eqref{T4h} implies that for any integer $N\geq 1$,
\be \ba \int_{\{|Tf|\geq 1\}}{|Tf(x)|^{\frac{nN}{n-\alpha}}}&=||F||_{q'N}^{q'N}\le {(2A_g)^{N}N!} \int_{\{|Tf|\geq 1\}}\exp\bigg[\frac{1}{2A_g}|T{f}(x)|^{\frac{n}{n-\alpha}}\bigg]dx\cr
&\le C{(2A_g)^{N}N!}||T{f}||_q^q\le C(2A_g)^{N}N^{N}||T{f}||_q^q .\label{T4i}\ea\ee
By \eqref{T4e}, we have $P_0,P_1\geq 1$, hence
\be \ba ||F||_{q'P_0}^{q'}&\le 2A_gC^{1/P_0}P_0||Tf||_q^{q/P_0}\le CP_0||Tf||_q^{q/P_0}\cr
||F||_{q'P_1}^{q'}&\le 2A_gC^{1/P_1}P_1||Tf||_q^{q/P_1}\le CP_1||Tf||_q^{q/P_1}.\ea\nt\ee
Let $a\in[0,1]$ be the number such that 
\be \frac{1}{P}=\frac a {P_0}+\frac{1-a}{P_1}.\nt\ee
By interpolation [Fol, Proposition 6.10] we have 
\be \ba ||F||_{q'P}^{q'}&\le ||F||_{q'P_0}^{q'a}||F||_{q'P_1}^{q'(1-a)}\le CP_0^aP_1^{1-a}||Tf||_q^{q(a/P_0)}||Tf||_q^{q(1-a)/P_1}\cr 
&= CP_0^aP_1^{1-a}||Tf||_q^{q/P}\le CP||Tf||_q^{q/P}.\ea\ee
Hence, since $p>1$
\be\ba I''&= \left(\int_{\{|Tf|\geq 1\}}{|Tf(x)|^{\frac{n}{n-\alpha}\frac{1}{p-1}}}dx\right)^{\frac{p-1}{p}}=||F||_{q'P}^{q'/p}\le C\bigg(\frac{1}{p-1}\bigg)^{1/p}||T{f}||_q^{q(\frac{p-1}{p})}\cr&\le \frac{C}{p-1}||T{f}||_q^{q(\frac{p-1}{p})}.
\ea\label{T4j}\ee
So combining \eqref{T4g} and \eqref{T4j}, we get
\be I'\cdot I''\le C\frac{1}{p-1}||T{f}||_q^{q}=C\frac{\tau}{p-1}\le C\frac{\tau}{1-(1-\tau)^{\frac{\a}{n-\a}}}\le C\label{T4k}\ee
where the second inequality is by \eqref{T4e}.
\section{Proof of Theorem \ref{T5} and Corollary \ref{corT5}}
\ \ \ In Theorem \ref{T1} we assume that the functions $f$ are compactly supported, with both $f$ and $Tf$ in the space $L^q(\R^n)$. We denote this space of functions by $D_0$:
\be D_0:=\{f\in L^q(\R^n)\ :\ \supp f\  \text{is compact and} \ Tf\in L^q(\R^n)\}.\nt\ee
In the following Theorem [FM2, Theorem 7] we see that $T$ has a smallest closed extension, which enables us to extend Theorem \ref{T1} to all functions in the domain of the extension. In particular, Theorem \ref{T5} is a consequence of the following Theorem:

 \bigskip
{\textup{\bf Theorem B ([FM2, Theorem 7])}.} \emph{If $K$ is a Riesz-like kernel, then the operator $T\ :\ D_0(T)\rightarrow L^q(\R^n)$ is closable, and its smallest closed extension (still denoted $T$) has domain
\be D(T)=\{f\in L^q(\R^n)\ :\ \exists \{f_k\}\subseteq D_0(T), \exists h\in L^q(\R^n) \ \text{with} \ f_k\xrightarrow[]{L^q} f, Tf_k\xrightarrow[]{L^q} h\}\label{T5a}\ee
and \be Tf=h.\nt\ee 
In the case of Riesz potential we have 
\be W^{\a,q}(\R^n)=\{I_\a f, f\in D(I_\a)\} \label{T5b}\ee
and the operator $(-\Delta)^{\frac{\a}{2}}$ is a bijection between $ W^{\a,q}(\R^n)$ and $D(I_\a)$, with inverse $c_\a I_\a.$}

 \bigskip
By using the above Theorem B and Fatou's lemma, we easily deduce that Theorem \ref{T1} is still valid for all functions $f$ in $D(T)$. Also \eqref{T5b} tells us that the Riesz potential for all functions $f$ in the extended domain $D(I_\a)$ is the space $W^{\a,q}(\R^n)$. Therefore by the fact that the inverse of $(-\Delta)^{\frac{\a}{2}}$ is $c_\a I_\a$, we have \eqref{5b}.\par
%
In the case of elliptic operator, by the formula \eqref{5f}, the kernel of the integral operator is homogeneous of order $\a-n$, therefore we have \eqref{5b}.\par
It is enough to assume $u\in C_c^\infty(\R^n)$ since $\a$ is an integer for the remaining cases. For $P=\nabla(-\Delta)^{\frac{\alpha-1}{2}}$ and $\alpha$ is an odd integer, since $u=c_{\a+1}I_{\a+1}(-\Delta)^{\frac{\alpha+1}{2}}u$, we can write
\be u(x)=\int_{\R^n} c_{\a+1}(n-\a-1)|x-y|^{\a-n-1}(x-y)\cdot f(y),\ \ \ f=\nabla(-\Delta)^{\frac{\alpha-1}{2}}u.\label{T5c}\ee
Clearly the kernel in the above formula satisfies our assumptions \eqref{A1}-\eqref{A4}, so \eqref{5b} follows.\par
For the proof of Corollary \ref{corT5}, it is clear that the inequality \eqref{cT5a} is a direct consequence of \eqref{5b} since $\Omega\subseteq \R^n$. 

 \bigskip
{\bf Proof of sharpness:}
To prove the sharpness, let $\psi_{\e,r}$ be the function as in the proof of sharpness (section \ref{sharpness}). If $P=(-\nabla)^{\frac{\a}{2}}$, consider the functions
\be u_{\e,r}=c_\a I_\a\psi_{\e,r}. \nt\ee
Similarly, for $P$ an elliptic operator, let $u_{\e,r}=g_P\ast \psi_{\e,r}$. \par
Lastly, we construct the extremal family of functions that proves sharpness for the case  $P=\nabla(-\Delta)^{\frac{\alpha-1}{2}}$ in Theorem \ref{T5}, as well as sharpness for Corollary \ref{corT5}. Note that in all these cases $\a$ is an integer. We use the same extremal functions as in Adams ([A1], see also [FM1], [FM2], [MS2]). Let $\varphi\in C^\infty([0,1])$ such that $\varphi^{(k)}(0)=0$ for $0\le k\le \a-1$, and $\varphi(1)=\varphi'(1)=1$, $\varphi^{(k)}(1)=0$ for $2\le k\le m-1$. Let $\e$ be small enough, define
\be v_\e(y)=\bc 0\qquad &\text{for}\ |y|\geq \frac{3}{4}\cr \varphi(\log\frac{1}{|y|})&\text{for}\ \frac{1}{2}\le|y|\le\frac{3}{4}\cr \log\frac{1}{|y|}&\text{for}\ 2\e\le|y|\le\frac{1}{2}\cr\log\frac{1}{\e}-\varphi(\log\frac{|y|}{\e}) &\text{for}\ \e\le|y|\le 2\e\cr\log\frac{1}{\e}&\text{for}\ |y|\le\e.\ec\nt\ee
Then we have that \be ||v_\e||_{\b q}\le C,\qquad ||\nabla^\a v_\e||_{q}^{q'}=\frac{\gamma(P)}{n}(\log\frac{1}{\e})^{q'-1}+O(1).\nt\ee
Let \be u_\e=\frac{v_\e}{||\nabla^\a v_\e||_{q}},\nt\ee
it is clear that 
\be ||\nabla u_\e||_q\le 1,\qquad ||u_\e||^{ q}_{q}\le C(\log\frac{1}{\e})^{-1},\label{mtshp1}\ee
and \be |u_\e|^{q'}\geq \gamma(P)^{-1}\log\frac{1}{\e^n},\qquad |y|\le\e.\nt\ee
For the sharpness of the exponential constant, we take $\theta >1$ and estimate
\be \ba \int_{\R^n}\frac{\exp_{\lceil\frac{n}{\alpha}-2\rceil}\big[\theta\gamma(P)|u_\e|^{\frac{n}{n-\alpha}}\big]}{1+|u_\e|^{\frac{ n}{n-\a}}}dy&\geq \int_{|y|\le\e}\frac{\exp\big[\theta\gamma(P)|u_\e|^{\frac{n}{n-\alpha}}\big]}{1+|u_\e|^{\frac{ n}{n-\a}}}dy\geq C\e^n\frac{\exp\bigg[\theta\log\dfrac{1}{\e^n}+C\bigg]}{1+C\log\dfrac{1}{\e}}\cr&= C\frac{\e^{(1-\theta)n}}{1+C\log\dfrac{1}{\e}}\to\infty \ea\nt\ee
as $\e\to0^+.$\par
For the sharpness of the power of the denominator, we take $\theta <1$ and get 
\be \ba \int_{\R^n}\frac{\exp_{\lceil\frac{n}{\alpha}-2\rceil}\big[\gamma(P)|u_\e|^{\frac{n}{n-\alpha}}\big]}{1+|u_\e|^{\frac{\theta n}{n-\a}}}dy&\geq \int_{|y|\le\e}\frac{\exp\big[\gamma(P)|u_\e|^{\frac{n}{n-\alpha}}\big]}{1+|u_\e|^{\frac{\theta n}{n-\a}}}dy\geq C\e^n\frac{\exp\bigg[\log\dfrac{1}{\e^n}+C\bigg]}{1+C(\log\dfrac{1}{\e})^{\theta}}\cr&\geq C(\log\frac{1}{\e})^{-\theta}. \ea\nt\ee
Hence by \eqref{mtshp1} we have that the quotient of the above integral over the norm $||u_\e||_{ q}^{ q}$ goes to infinity as $\e\to0^+$, so the sharpness follows.\par



\vfill\eject
\appendix
\ni{\LARGE \bf Appendices}\bigskip
\section{Proof of Lemma \ref{l}  (Improved O'Neil Lemma) }\label{lemma1}

\bigskip

For simplicity assume that the hypothesis on $f_x$ are true for all $x\in N$; the proof works just the same if $N$ is replaced by $N\setminus F$, some measurable $F$ with $\nu(F)=0$. First let us recall the following special case of a weak-type estimate due to Adams [A3, Lemma1]: if $k$ is nonnegative and satisfies \eqref{l1} then for any nonnegative $\phi\in L^1$ and any $s>0$
\be s\,\nu\big(\{x\in N:T\phi(x)>s\}\big)^{\frac1{\beta}}\le \b' D\|\phi\|_1.\label{wt1}\ee
In this special case the above estimate can be proved easily: take any $X\subseteq N$ with finite measure and  let $E_s=\{x\in X: T\phi(x)>s\}.$ Then we have 
\be s\,\nu(E_s)\le \int_M\phi(y)d\mu(y)\int_{E_s}k(x,y)d\nu(x)\le \|\phi\|_1\int_0^{\nu(E_s)} k_1^*(t)dt\le\beta' D\|\phi\|_1\nu(E_s)^{\frac1{\beta'}}\nt\ee
which, using $\sigma-$ finiteness of $(N,\nu)$,  gives \eqref{wt1}.

Returning to Lemma 1, we can apply \eqref{wt1} to the kernel $|k(x,y)|$ and  $\phi= \ol f$, and  get that $T\ol f(x)<\infty$ for  $\nu-$a.e. $x$, and hence $T'f(x)$ is well-defined and finite for  $\nu-$a.e. $x$.
 
Without loss of generality we can now assume that both $f$ and $k$ are nonnegative.

Fix $t,\tau>0$, pick a sequence $\{s_n\}_{-\infty}^{\infty}$ such that $s_0=\ol{f}^*(\tau),\ s_n< s_{n+1},\ s_n\rightarrow \infty$ as $n\rightarrow \infty$, and $s_n\rightarrow 0$ as $n\rightarrow -\infty$. Define for each $x\in N$, $y\in M$
\be f_{n}(y)=
\bc
0\ \ \ &\text{if} \ f(x,y)\leq s_{n-1}\\
f(x,y)-s_{n-1}\ &\text{if} \ s_{n-1}<f(x,y)\leq s_{n}\\
s_n-s_{n-1} \ &\text{if} \ s_n<f(x,y)
\ec\nt\ee\be\ol{f}_{n}(y)=
\bc
0\ \ \ &\text{if}\ \ol{f}(y)\leq s_{n-1}\\
\ol{f}(y)-s_{n-1}\ &\text{if} \ s_{n-1}<\ol{f}(y)\leq s_{n}\\
s_n-s_{n-1} \ &\text{if} \ s_n<\ol{f}(y)
\ec
\nt\ee

so that 
\be f=\sum_{-\infty}^0 {f}_{n}+\sum_1^{\infty} {f}_{n}:=g_1+g_2,\qquad \ol{f}=\sum_{-\infty}^0 \ol{f}_{n}+\sum_1^{\infty} \ol{f}_{n}:=\ol g_1+\ol g_2 \label{l8}\ee\def\supp{\text {supp}}
Denoting $f_{x,n}(y)=f_n(x,y),\; g_{x,j}=g_j(x,y)$, it is clear that  from (2.3) we have, for all $x\in N$,  
\be g_{x,2}\leq \ol g_2,\ \ \mu-a.e.\label{l9}\ee
Letting \be T'g_2(x)=Tg_{x,2}(x)=\int_M k(x,y)g_{x,2}(y)d\mu(y) \nt\ee
then we have, for all $x\in N$,  \be T'g_2(x)\leq T\ol g_2(x)\label{l10}\ee
and therefore, by the subadditivity of $(\cdot)^{**}$
 \be (T'g_2)^{**}(t)\leq (T\ol g_2)^{**}(t)\leq \sum_1^{\infty} (T\ol{f}_n)^{**}(t).\ee
 Now note that the weak type estimate \eqref{wt1} is equivalent to $(T\phi)^*(t)\le \beta' D t^{-\frac1{\beta}}\|\phi\|_1$. Applying this fact to $\phi=\ol f_n$ we get
 \be\ba (T\ol f_n)^{**}(t)\leq Ct^{-\frac{1}{\beta}}||\ol f_n||_1&\le Ct^{-\frac{1}{\beta}}(s_n-s_{n-1})\mu(\supp\ol f_n)\cr &=Ct^{-\frac{1}{\beta}}(s_n-s_{n-1})m_{\ol f}(s_{n-1})\ea\ee
 and 
 \be (T'g_2)^{**}(t)\leq C t^{-\frac{1}{\beta}}\sum_1^\infty(s_n-s_{n-1})m_{\ol f}(s_{n-1}).\ee
 Taking the inf over all such $\{s_n\}$ in the above Riemann sum we get
 \be 
 \ba
 &(T'g_2)^{**}(t)\le (T\ol g_2)^{**}(t)\leq Ct^{-\frac{1}{\beta}}\int_{\ol{f}^*(\tau)}^\infty m_{\ol{f}}(s)ds=Ct^{-\frac{1}{\beta}}\int_{\ol{f}^*(\tau)}^\infty m_{\ol{f}^*}(s)ds\cr
 &=Ct^{-\frac{1}{\beta}}\int_{\ol{f}^*(\tau)}^\infty ds\int_{\{\ol f^*(u)>s\}} du=Ct^{-\frac{1}{\beta}}\int_0^\infty du\int_{\{\ol f^*(\tau)\le s<\ol f^*(u)\}} ds\cr
 &=Ct^{-\frac{1}{\beta}}\int_0^\tau du\int_{\{\ol f^*(\tau)\le s<\ol f^*(u)\}} ds\leq Ct^{-\frac{1}{\beta}}\int_0^\tau \ol{f}^*(u)du.
 \ea
 \label{tg2}\ee

To deal with $g_1$, we fix $ x\in N$ and write
\be\ba T'g_1(x)&:=Tg_{x,1}(x)\leq\sum_{-\infty}^0 Tf_{x,n}(x)\leq \sum_{-\infty}^0(s_n-s_{n-1})\int_{\supp f_{x,n}}k(x,y)d\mu(y)\cr& \leq \sum_{-\infty}^0(s_n-s_{n-1})\int_0^{m_{f_{x}^*}(s_{n-1})}k_1^*(x,u)du.
\ea
\label{l13}\ee
Passing once again to the limit in the partition $\{s_n\}$ we obtain

\be T'g_{1}(x)\leq\int_0^{\ol f^*(\tau)}ds\int_0^{m_{f_{x}^*}(s)}k_1^*(x,u)du= \int_0^\infty \!k_1^*(x,u)\,\varphi(u)du.
 \ee

where \be \varphi(u)=\big|\big\{s: 0<s<\ol f^*(\tau),\,u<m(f_x^*,s)\big\}\big|\le\min\big\{\ol f^*(\tau),\,  f_x^*(u)\big\}\label{fx}\ee
To get the above estimate note that if $s\ge f_x^*(u)$ then $m_{f_x^*}(s)\le m_{f_x^*}(f_x^*(u))\le u.$

Hence we have
\be\ba T'g_{1}(x)&\leq \bigg(\int_0^\tau\!\!+\!\!\int_\tau^\infty\bigg) k_1^*(x,u)\,\varphi(u)du\le \ol f^*(\tau) \int_0^\tau k_1^*(x,u)du+\int_\tau^\infty k_1^*(x,u)f_x^*(u)du\cr&\le C\tau^{-\frac1{\beta}}\int_0^\tau \ol f^*(u)du+\int_\tau^\infty k_1^*(x,u)f_x^*(u)du.\ea\label{tg1}\ee 

\ni Putting together \eqref{tg2} and \eqref{tg1}, from $T'f(x)=T'g_1(x)+T' g_2(x)$ we obtain
\be\ba (T'f)^{**}(t)&\le (T\ol g_2)^{**}(t)+(T'g_1)^{**}(t)\le Ct^{-\frac{1}{\beta}}\int_0^\tau \ol{f}^*(u)du+\|T'g_1\|_\infty\cr&
\le C\max\big\{\tau^{-\frac1{\beta}},t^{-\frac{1}{\beta}}\big\}\int_0^\tau \ol{f}^*(u)du+\supess_{x\in N}\int_\tau^\infty k_1^*(x,u)f_x^*(u)du.\ea\ee

 \hspace*{\fill} \QEDopen\\

\ni{\bf Remark 4.} The computations involved in  \eqref{tg2} and \eqref{tg1} are more direct than those in the existing literature for the case  $f(x,y)=f(y)=\ol f(y)$ (see for example [FM1], [FM3]). We make use of Fubini's theorem, whereas existing proofs, based on O'Neil's original argument, make use of  integrals with respect to $df^*$ and integration by parts (see  [Z], proof of Lemma 1.8.8, for a detailed justification of those steps). The simplifications in \eqref{tg2}-\eqref{tg1} were suggested to me by Luigi Fontana and  Carlo Morpurgo, to whom I am grateful.


\section{Proof of Lemma \ref{lb}}\label{lemma6}

\bigskip
\ni{\bf Proof of \eqref{lb1}:} Using \eqref{A3} we get
  \be \ba \avgdJ |Tf_\tau(x)|dx
  &\leq \frac{C}{|D_{J+1}\setminus D_J|}\int _{\R^n}|f_\tau(y)|\int_{D_{J+1}\setminus D_J}|x-y|^{\a-n}dxdy\cr&\leq  \frac{C}{r_J^n}\int _{\R^n}|f_\tau(y)|r_J^\a dy. \ea\label{lb2}\ee
  Here the second inequality above is by the straightforward computation:
  \be \ba &\int_{D_{J+1}\setminus D_J}|x-y|^{\a-n}dx\cr
  &=\int_{\{|x-y|\leq r_{J}\}\cap (D_{J+1}\setminus D_J) }|x-y|^{\a-n}dx+\int_{\{|x-y|> r_{J}\}\cap (D_{J+1}\setminus D_J)}|x-y|^{\a-n}dx\cr
  &\le \int_{\{|x|\leq r_{J}\}}|x|^{\a-n}dx+\int_{D_{J+1}\setminus D_J}r_J^{\a-n}dx\le Cr_J^\a+Cr_J^{\a-n}r_J^n=Cr_J^\a. \ea\label{lb222}\ee
  Recall that $|F_\tau|=\tau=|D_0|$, we have
  \be \ba &\frac{C}{r_J^n}\int _{\R^n}|f_\tau(y)|r_J^\a dy=Cr_J^{\a-n}\int _{\R^n}|f_\tau(y)|dy=Cr_J^{\a-n}\int _{\R^n}|f|\chi^{}_{F_\tau}dy\cr
  &\le Cr_J^{\a-n}|F_\tau|^{1/q'} ||f||_q\leq Cr_J^{\a-n}r_0^{n-\a}=C_1\left(\frac{r_0}{r_J}\right)^{n-\a}=C_1\left(\frac{1}{e^{q-1}}\right)^J.\ea\label{lb3}\ee
  \\
   \ni{\bf Proof of \eqref{lc1}:} First write 
 \be {T}S_{J+2}f_\tau'(x)-{T}S_{J+1}f_\tau'(x_1)={T}S_{J+2}f_\tau'(x)-{T}S_{J+2}f_\tau'(x_1)-{T}\big(S_{J+1}-S_{J+2})f_\tau'(x_1)\nt\ee
 so \be \ba &| {T}S_{J+2}f_\tau'(x)- {T}S_{J+1}f_\tau'(x_1)|\cr
 &\le|{T}S_{J+2}f_\tau'(x)-{T}S_{J+2}f_\tau'(x_1)|+|{T}\big(S_{J+1}-S_{J+2})f_\tau'(x_1)|.
 \ea\label{lc2}\ee
 Arguing as in \eqref{la3}, we get 
 \be \ba |{T}\big(S_{J+1}-S_{J+2})f_\tau'(x_1)|&=\bigg|\int_{D_{J+2}\setminus D_{J+1}}K(x_1-y)f_\tau'(y)dy\bigg|\cr
 &\leq C\bigg|\int_{D_{J+2}\setminus D_{J+1}}|x_1-y|^{\a-n}f_\tau'(y)dy\bigg|\cr&\leq Cr_{J+1}^{\a-n}\int_{D_{J+2}\setminus D_{J+1}}|f_\tau'(y)|dy\le Cr_{J+1}^{\a-n}r_{J+1}^{n-\a}\a_{J+1}\le C\b_{J+1}.
 \ea\label{lc3}\ee
By the regularity assumption \eqref{A3}, since $x_1\in D_0$, we have for $x\in D_{J+1}\setminus D_{J}$ and $y\in D^c_{J+2}$
 \be |K(x-y)-K(x_1-y)|\leq C|x-x_1|(e^{q/n})^{n+1-\a}|x_1-y|^{\a-n-1}.\nt\ee
 Hence 
 \be \ba  |{T}S_{J+2}f_\tau'(x)-{T}S_{J+2}f_\tau'(x_1)|
& =\bigg|\int_{D^c_{J+2}}(K(x-y)-K(x_1-y))f_\tau'(y)dy\bigg|\cr
&\leq C|x-x_1|\int_{D^c_{J+2}}|f_\tau'(y)||x_1-y|^{\a-n-1}dy\cr
&\leq C|x-x_1|\b_{J+2}\left(\int_{D^c_{J+2}}|x_1-y|^{-n-\frac{n}{n-\a}}dy \right)^{\frac{n-\a}{n}}\cr
&\leq Cr_{J+1}\b_{J+2}\frac{C}{r_{J+2}}\leq C\b_{J+2}\le C\b_{J+1}.\ea\label{lc4}\ee
So we have \eqref{lc1} by \eqref{lc2}-\eqref{lc4}.

 \bigskip
\ni {\bf Proof of \eqref{ld1} and \eqref{ld2}:} 
 Let $j\in\{0,1,...,J\}$, then
\be \ba\bigg|\avgdJ T(S_j-S_{j+1})f_\tau'(x)dx\bigg|&=\bigg|\avgdJ \dj K(x-y)f_\tau'(y)dydx\bigg|\cr
&\le C\avgdJ \dj |x-y|^{\a-n}|f_\tau'(y)|dydx\cr
&=C\dj |f_\tau'(y)|\avgdJ |x-y|^{\a-n}dxdy\cr
&\le Cr_J^{\a-n}\dj |f_\tau'(y)|dy\le Cr_J^{\a-n} r_j^{n-\a}\a_j\cr&=Ce^{\frac{q}{n}(j-J)}\a_j \ea
\label{ld3}\ee
where the second inequality is by \eqref{lb222}.
Therefore,
\be \ba \bigg|\avgdJ{T}(S_0-S_{J+2})f_\tau'(x)dx\bigg|&=\bigg|\avgdJ \left(\sum_{j=0}^{J+1}T(S_j-S_{j+1})f_\tau'(x)\right)dx\bigg|\cr
&\le \sum_{j=0}^{J+1}\bigg|\avgdJ T(S_j-S_{j+1})f_\tau'(x)dx\bigg|\cr
&\le  C\sum_{j=0}^{J+1} e^{\frac{q}{n}(j-J)}\a_j\le C\ab_{J+1}.
\ea
\label{ld4}\ee
So we have \eqref{ld1}. For \eqref{ld2}, by calculations similar to those in \eqref{ld3}, we have
\be \ba \bigg|\avgdJ T(f_\tau'\chi_{D_0})(x)dx\bigg|
&=\bigg|\avgdJ \int_{D_0} K(x-y)f_\tau'(y)dydx\bigg|\cr
&\le C\avgdJ \int_{D_0} |x-y|^{\a-n}|f_\tau'(y)|dydx\cr
&=C\int_{D_0}|f_\tau'(y)|\avgdj |x-y|^{\a-n}dxdy\cr
&\le Cr_J^{\a-n}\int_{D_0}|f_\tau'(y)|dy\le Cr_J^{\a-n} r_0^{n-\a}\a_{-1}\cr
&\le C\a_{-1}. \ea
\label{ld5}\ee
\hfill\QEDopen
\vskip 1in
\ni{\bf Acknowledgments.} The results presented in this paper are part of the author's Ph.D. dissertation [Q] at University of Missouri, Columbia. The author is grateful to her advisor Carlo Morpurgo for his advice and useful suggestions.

\bigskip\bigskip
\ni Liuyu Qin

\ni Department of Mathematics and statistics

\ni Hunan University of Finance and Economics

\ni Changsha, Hunan

\ni China

\smallskip\ni \text{Liuyu\_Qin@outlook.com}

\end{document}